\documentclass[10pt]{elsarticle}
\usepackage{lineno,hyperref}
\usepackage{booktabs}
\usepackage{multirow}
\usepackage{amsthm}
\usepackage{graphicx}
\usepackage{subfigure}
\usepackage{float}
\usepackage{bm}
\usepackage{lipsum}
\usepackage{listings} 
\usepackage{color}
\usepackage{amsfonts,enumerate,amsmath,amssymb}
\def\qed{\strut\hfill $\Box$}
\newtheorem{thm}{Theorem}[section]

\newtheorem{prop}[thm]{Proposition}
\newtheorem{lem}[thm]{Lemma}

\newtheorem{rem}[thm]{Remark}
\newtheorem{defn}[thm]{Definition}

\newcommand{\thmref}[1]{Theorem~{\rm \ref{#1}}}
\newcommand{\lemref}[1]{Lemma~{\rm \ref{#1}}}
\newcommand{\propref}[1]{Proposition~{\rm \ref{#1}}}
\newcommand{\defref}[1]{Definition~{\rm \ref{#1}}}
\newcommand{\remref}[1]{Remark~{\rm \ref{#1}}}

\def\para#1{\vskip .4\baselineskip\noindent{\bf #1}}
\newcommand{\vertiii}[1]{{\left\vert\kern-0.25ex\left\vert\kern-0.25ex\left\vert #1 
		\right\vert\kern-0.25ex\right\vert\kern-0.25ex\right\vert}}
\def\para#1{\vskip .4\baselineskip\noindent{\bf #1}}
\numberwithin{equation}{section}
\begin{document}
	\begin{frontmatter}

		\title{Large deviation principle for slow-fast rough differential equations via controlled rough paths}
		
		\author[mymainaddress]{Xiaoyu Yang}
		\ead{yangxiaoyu@yahoo.com}

		\author[mysecondaddress]{Yong Xu\corref{mycorrespondingauthor}}\cortext[mycorrespondingauthor]{Corresponding author}\ead{hsux3@nwpu.edu.cn}

		\address[mymainaddress]{Graduate School of Information Science and Technology,  Osaka University, Osaka, 5650871, Japan}
		\address[mysecondaddress]{School of Mathematics and Statistics, Northwestern Polytechnical University, Xi'an, 710072, China}

		\begin{abstract}
We prove a large deviation principle for the slow-fast rough differential equations under the controlled rough path framework. The driver rough paths are lifted from the mixed fractional Brownian motion with Hurst parameter $H\in (1/3,1/2)$. Our approach is based on the continuity of the solution mapping and the variational framework for mixed fractional Brownian motion. By utilizing the variational representation, our problem is transformed into a qualitative property of the controlled system. In particular, the fast rough differential equation coincides with It\^o SDE almost surely, which possesses a unique invariant probability measure with frozen slow component. We then demonstrate the weak convergence of the controlled slow component by averaging with respect to the invariant measure of the fast equation and exploiting the continuity of the solution mapping. 
			
			\vskip 0.08in			
				\noindent{\bf AMS Math Classification.}
			{60F10, 60G15, 60H10}.
			\vskip 0.08in			
		\end{abstract}		
	\end{frontmatter}

	\section{Introduction}\label{sec-1}
	The topic of this paper is to studying the  slow-fast rough differential equation (RDE) in time interval $[0,T]$ under the controlled rough path (RP) framework as follows:
	\begin{eqnarray}\label{1}
	\left
	\{
	\begin{array}{ll}
	X^{\varepsilon, \delta}_t =X_0+\int_{0}^{t} f_{1}(X^{\varepsilon, \delta}_s, Y^{\varepsilon, \delta}_s)ds + \int_{0}^{t}\sqrt \varepsilon  \sigma_{1}( X^{\varepsilon, \delta}_s)dB^H_{s},\\
	Y^{\varepsilon, \delta}_t =Y_0+\frac{1}{\delta} \int_{0}^{t}f_{2}( X^{\varepsilon, \delta}_s, Y^{\varepsilon, \delta}_s)ds + \frac{1}{{\sqrt \delta}}\int_{0}^{t}\sigma_2( X^{\varepsilon, \delta}_s, Y^{\varepsilon, \delta}_s)dW_{s}.
	\end{array}
	\right.
	\end{eqnarray}
	Here, 	the RP $(B,W)$  is lifted from the mixed fractional Brownian motion (FBM) $(b^H,w)$ with Hurst parameter $H\in (\frac{1}{3},\frac{1}{2})$.  Two small parameters $\varepsilon$ and $\delta$ satisfies the condition that $0<\delta=o(\varepsilon)<\varepsilon\leq 1$.  ${X^{\varepsilon,\delta}}$ is  the slow component and ${Y^{\varepsilon,\delta}}$ is the fast component with the  (arbitrary but deterministic) initial data $(X^{\varepsilon, \delta}_0, Y^{\varepsilon, \delta}_0)=(X_0, Y_0)\in \mathbb{R}^{{m}}\times \mathbb{R}^{{n}}$. The coefficients $f_{1}:\mathbb{R}^{{m}}\times \mathbb{R}^{{n}}\to\mathbb{R}^{{m}}$, $f_2:\mathbb{R}^{{m}}\times \mathbb{R}^{{n}}\to\mathbb{R}^{{n}}$, {$\sigma_{1}: \mathbb{R}^{{m}}\to \mathbb{R}^{{m}\times{d}}$ and $\sigma_{2}: \mathbb{R}^{{m}}\times \mathbb{R}^{{n}}\to \mathbb{R}^{{n}\times{e}} $} are nonlinear regular enough functions, which assumed to satisfy some suitable conditions   in section 3. 	Such a slow-fast model has been applied in many real world fields, for example,   typical examples could be found in  climate-weather  (see  \cite{2000Kiefer}), biological field and so on \cite{Krupa2008Mixed}.   The dynamical behavior for slow-fast model is a active research area, see for instance, the monographs  \cite{2008Pavliotis} and references \cite{2021Bourguin,2020Hairer,2023Rockner} therein for a comprehensive overview. 

As a generalization of the standard Wiener process $(H=1/2)$, the FBM is self-similar and possesses long-range dependence, which has become  widely popular for applications \cite{2008Biagini,2002Rascanu,1999Decreusefond}. Its Hurst parameter $H$ could depict   the roughness of the sample paths, with a lower value leading to a rougher motion \cite{2008Mishura}. Especially,   the case of $H<1/2$ seems  rather troublesome to be handled with the conventional stochastic techniques. To get over the hump that is caused by rougher sample paths for $H<1/2$, our model  is within the RP setting.	The so-called RP theory does not require martingale theory, Markovian property, or filtration theory.  This also determines the de-randomisation  when being applied in the stochastic situation, so it can provide a new prescription  to FBM problems.   The RP theory was originally proposed by Terry Lyons in 1998 \cite{1998Lyons,2002Lyons} and has sparked  tremendous  interest from the fields of probability \cite{2013Inahama,2015Inahama} and applied mathematics \cite{2010Friz,2019Liu} after 2010.
Briefly, the main idea of RP theory states that it not only considers the path itself, but also considers the iterative integral of the path, so that the continuity of the solution mapping could be ensured. This continuity property of the solution mapping is the core of RP theory.  Until now, there have been three formalisms to RP theory \cite{2020Friz,2002Lyons,2010Friz} and  we adopt that one of them, which is so-called controlled RP theory  \cite{2020Friz}.
	 By resorting to the controlled RP framework, the slow-fast RDE (\ref{1}) under suitable conditions admits a unique (pathwise) solution $(X^{\varepsilon, \delta},Y^{\varepsilon,\delta})\in \mathcal{C}^{\beta-\operatorname{hld}}([0,T],\mathbb{R}^{m}) \times \mathcal{C}\left([0,T], \mathbb{R}^n\right)$ with $\beta\in (0,H)$, which will be precisely stated in Section 3.
		Here, $\mathcal{C}^{\beta-\operatorname{hld}}([0,T],\mathbb{R}^{m})$ and
			$\mathcal{C}\left([0,T], \mathbb{R}^n\right)$ are the $\beta$-H\"older continuous 
			path space and the continuous path space, respectively.	  	  
 
	In accordance with the averaging principle, as $\delta \to 0$, $X^{\varepsilon, \delta}$ 	is well approximated by an effective dynamics $\bar X$ which is defined as following, 
\begin{equation}\label{1-4}
	\left\{\begin{array}{l}
	d \bar{X}_t=\bar{f}_1\left(\bar{X}_t\right) d t \\
	\bar{X}_0=X_0 \in \mathbb{R}^m,
	\end{array}\right.
	\end{equation}
	with $\bar{f}_1(x)=\int_{\mathbb{R}^{n}} f_1(x, y) d \mu_x(y)$. Here,  $\mu_x$ is a unique invariant probability measure of the fast component with the ``frozen"-$x$. The  precise proof is a small extension of \cite[Theorem 2.1]{2022Inahama}. 
	
	However, the small parameter $\delta$ can not be zero and when it is small enough, the trajectory of the slow component would  stay in a small neighborhood of $\bar X$.
The Large Deviation Principle (LDP) could describe the extent to which the slow component deviates from the average component exponentially, which is more accurate. 	
	As a result, the main objective of  this work is to prove a LDP for the slow component $X^{\varepsilon,\delta}$ of the above RDE \eqref{1}.  
	The family
	${X^{\varepsilon,\delta}}\in \mathcal{C}^{\beta-\operatorname{hld}}([0,T],\mathbb{R}^{m})$ is called to satisfy a LDP on $\mathcal{C}^{\beta-\operatorname{hld}}([0,T],\mathbb{R}^{m})$ ($0 <\beta <H$) with a good
	rate function $I: \mathcal{C}^{\beta-\operatorname{hld}}([0,T],\mathbb{R}^{m})\rightarrow [0, \infty]$  if the following two conditions hold:
	\begin{itemize}
		\item For each closed subset $F$ of $\mathcal{C}^{\beta-\operatorname{hld}}([0,T],\mathbb{R}^{m})$,
		$$\limsup _{\varepsilon \rightarrow 0} \varepsilon \log \mathbb{P}\big(X^{\varepsilon,\delta}\in F\big) \leqslant-\inf _{x \in F} I(x).$$
		\item For each open subset $G$ of $\mathcal{C}^{\beta-\operatorname{hld}}([0,T],\mathbb{R}^{m})$,
		$$\liminf _{\varepsilon \rightarrow 0} \varepsilon \log \mathbb{P}\big(X^{\varepsilon,\delta}\in G\big) \geqslant-\inf _{x \in G} I(x).$$
	\end{itemize}
	This will be stated in our main result (Theorem \ref{thm}) and  the definition of $I$ will also be given there. 	     

The LDP for stochastic dynamical systems was pioneered by Freidlin and Wentzell \cite{1984Random}, which has inspired much of the subsequent substantial development \cite{2004Cerrai,2012Dupuis,1994Peszat,2006Sritharan}. Up to date, there have been several different approaches to studying LDP for the stochastic slow-fast system, such as the weak convergence method \cite{2011Variational,2011Dupuis,2020Large}, the PDE theory \cite{BCS}, nonlinear semigroups and viscosity solution theory \cite{FFK,FK}. It is remarkable that the weak convergence method, which is founded on the variational representation for the non-negative functional of BM \cite{1998Dupuis}, has been extensively utilised in the LDP of the slow-fast systems with BM. As well as this, the weak convergence method is powerful for solving LDP problems in FBM situations \cite{2020Budhiraja, 2023Inahama} with $H>1/2$.
	
	Nevertheless, it is a priori not clear if the  LDP  for slow-fast RDE \eqref{1} holds and  the aforementioned  methods are not sufficient to answer this question.  For the single-time scale RDE, the RP theory is proven  efficient in the LDP problems  by using the exponentially good approximations of Gaussian processes \cite{2007Friz, 2002Ledoux, 2006Millet}. However, due to hinging on the fast equation, this exponentially good approximation method is invalidated in our slow-fast case. In response to this challenge, new approach has to be developed. Our work is  to  adopt the variational framework to solve the LDP for the slow-fast RDE. The technical core of the proof is    the continuity of the solution mapping and the weak convergence method, which is  based on the  variational representation of mixed FBM. Before stating our outline of the proof,  we firstly give two important results.  Thanks to  \cite[Proposition 4.7]{2022Inahama}, it is pointed out that for each $0<\delta, \varepsilon\leq 1$, $Y^{\varepsilon,\delta}$ coincides with the It\^o SDE almost surely and it possesses a unique invariant probability measure with frozen slow component.  The second result is that  the translation of mixed FBM in the direction of Cameron-Martin components can be lifted to RP, which will be proved in section 2.  Then, we give the outline of our proof. Firstly, based on the  variational representation formula for a standard BM  \cite{BP_book},   the variational representation formula for  mixed  RP is given. Then, the LDP problem  could be transformed into weak convergence of the  controlled slow RDE. It is a key ingredient  in  the  weak convergence  to  average out the controlled fast component.    Then, we  show that the controlled fast component could be replaced by the fast component without controlled term in the limit by the condition that $\delta=o(\varepsilon)$. Finally, we derive the weak convergence of the controlled slow component by exploiting the exponential ergodicity of the auxiliary fast component without control, continuity of the solution mapping, the continuous mapping theorem, and so on.

	We now give the outline of this paper. In the  Section 2, we introduce some notation and   preliminaries.  In the  Section 3, we give  assumptions and the statement of our main result. Section 4 is devoted to  prior estimates.  In Section 5,  the proof of our main result is achieved. 
	Throughout this paper, $c$, $C$, $c_1$, $C_1$, $\cdots$ denote certain positive constants that may vary from line to line. {$\mathbb{N}=\{1,2, \ldots\}$} and time horizon $T >0$. 
	
	\section{Notations and Preliminaries}\label{sec-2}
	\textbf{Notations}
	Firstly, we introduce the  notations which will be used throughout our paper.
	Let $[a,b]\subset[0,T]$ and $\Delta_{[a,b]}:=\{(s,t)\in \mathbb{R}^2|a\leq s\leq t\leq b\}$. We write $\Delta_{T}$ simply when $[a,b]=[0,T]$.  Denote $\nabla$ be the  standard gradient on a Euclidean space. {Throughout this section, $\mathcal{V}$ and $\mathcal{W}$ are Euclidean spaces.}  
		\begin{itemize}
		\item 	 \textbf{(Continuous space)} Denote $\mathcal{C}([a,b],\mathcal{V})$ by the space of continuous functions $\varphi:[a,b]\to \mathcal{V}$  with the norm $\|\varphi\|_{\infty}=\sup _{t \in[a, b]}|\varphi_t|<\infty$, which is a Banach space. The set of  continuous functions starts from $0$ is denoted by $\mathcal{C}_0([a,b],\mathcal{V})$.
		\item 	   \textbf{(H\"older continuous space and variation space)}
		For $\eta \in(0,1]$,  denote $\mathcal{C}^{\eta-\operatorname{hld}}([a,b],\mathcal{V})$ by the space of $\eta$-H\"older continuous functions $\varphi:[a,b]\to \mathcal{V}$, equipped with the semi-norm
		$$
		\|\varphi\|_{\eta-\operatorname{hld},[a,b]}:=\sup _{a \leq s<t \leq b} \frac{|\varphi_t-\varphi_s|}{(t-s)^\eta}<\infty.
		$$ 
		The Banach norm in  $\mathcal{C}^{\eta-\operatorname{hld}}([a,b],\mathcal{V})$ is $ |\varphi_a|_{\mathcal{V}}+
		\|\varphi\|_{\eta-\operatorname{hld},[a,b]}$.
		
		For $1\leq p<\infty $,  denote $\mathcal{C}^{p-\operatorname{var}}\left([a,b],\mathcal{V}\right)=\left\{\varphi \in \mathcal{C}\left([a,b],\mathcal{V}\right):\|\varphi\|_{p \text {-var }}<\infty\right\}$ where $\|\varphi\|_{p \text {-var }}$ is   the usual $p$-variation semi-norm. The set of $\eta$-H\"older continuous functions starts from $0$ is denoted by $\mathcal{C}_0^{\eta-\operatorname{hld}}([a,b],\mathcal{V})$.	The space $\mathcal{C}_0^{p \text{-var}}([a,b],\mathcal{V})$ is defined in a similar way.
		
		For a continuous map $\psi:\Delta_{[a,b]}\to \mathcal{V}$, we set 
		$$		\|\psi\|_{\eta-\operatorname{hld},[a,b]}:=\sup _{a \leq s<t \leq b} \frac{|\psi_t-\psi_s|}{(t-s)^\eta}.$$
		We denote the set of above such $\psi$ of $\|\psi\|_{\eta-\operatorname{hld},[a,b]}<\infty$  by $\mathcal{C}_2^{\eta-\operatorname{hld}}([a,b],\mathcal{V})$. It is a Banach space equipped with the norm $\|\psi\|_{\eta-\operatorname{hld},[a,b]}$.

\item \textbf{(Besov space)} For $\phi:[a,b] \rightarrow \mathcal{V}$   and  $\delta \in(0,1) $  and  $p \in(1, \infty)$,  we define the Besov space $W^{\delta, p}\left([a,b],\mathcal{V}\right)$ equipped with the following norm:
\begin{equation}\label{2-2}
\|\phi\|_{W^{\delta, p}}=\|\phi\|_{L^p}+\left(\iint_{[a,b]^2} \frac{\left|\phi_t-\phi_s\right|^p}{|t-s|^{1+\delta p}} d s d t\right)^{1 / p} <\infty.
\end{equation}
Moreover, when $\eta'=\delta-1/p>0$, we have the continuous imbedding that $W^{\delta, p}([a,b],\mathcal{V}) \subset \mathcal{C}^{\eta'-\operatorname{hld}}([a,b],\mathcal{V})$ \cite[Theorem 2]{2006Friz}.
		
	\item	\textbf{{($C^k$ norm and $C^k_b$ norm)}}
	Let $U\subset \mathcal{V}$ be an open set. For $k\in \mathbb{N}$, denote $C^k(U,\mathcal{W})$ by the set of $C^k$-functions from $U$ to $\mathcal{W}$. $C_b^k(U,\mathcal{W})$ stands  the set of $C^k$-bounded functions  whose  derivatives up to order $k-$ are also bounded. The space $C_b^k(U,\mathcal{W})$ is a Banach space equipped with the norm $\|\varphi\|_{C_b^k}:=\sum_{i=0}^{k}\|\nabla^i\varphi\|_\infty<\infty$.
	\item $L(\mathcal{W},\mathcal{V})$ denotes the set of bounded linear maps     from $\mathcal{W}$ to $\mathcal{V}$. We set $L(\mathcal{V}, L(\mathcal{V}, \mathcal{W})) \cong L^{(2)}(\mathcal{V} \times \mathcal{V}, \mathcal{W}) \cong L(\mathcal{V} \otimes \mathcal{V}, \mathcal{W})$ where $L^{(2)}(\mathcal{V} \times \mathcal{V}, \mathcal{W})$ is the vector space of bounded bilinear maps from $\mathcal{V} \times \mathcal{V}$ to $\mathcal{W}$.
		\end{itemize}  
	\subsection{Mixed fractional Brownian motion}\label{sec-2-1}
	This subsection features a brief overview of the mixed FBM of Hurst prameter $H$, and  only focuses on the case that  $H\in(1/3,1/2)$.
		
		Consider the $\mathbb{R}^{d}$-valued continuous stochastic process $(b^H_t)_{t\in[0,T]}$ starting from 0 as following:
	$$b^H_t=(b_t^{H,1},b_t^{H,2},\cdots,b_t^{H,d}).$$
	The above  $(b^H_t)_{t\in[0,T]}$ is said to be a FBM  if it is  a centred Gaussian process, satisfying that
	$$\mathbb{E}\big[b_{t}^{H} b_{s}^{H}\big]=\frac{1}{2}\left[t^{2 H}+s^{2 H}-|t-s|^{2 H}\right]\times I_{d}, \quad(0\leq s\leq t \leq T),$$
	where $I_{d}$ stands  the identity matrix in $\mathbb{R}^{d\times d}$. 
	Then, it is easy to see that  
	$$\mathbb{E}\big[(b_{t}^{H}-b_{s}^{H})^{2}\big]=|t-s|^{2 H} \times I_{d},\quad (0\leq s\leq t \leq T).$$
From the
Kolmogorov's continuity criterion, 
the trajectories of $b^H$  is of $H'$-H\"older continuous ($H'\in(0,H)$) and $\lfloor 1 / H \rfloor<{p}<\lfloor 1 / H \rfloor+1$-variation almost surely. 
	The reproducing kernel Hilbert space of the FBM $b^H$ is denoted by $\mathcal{H}^{H,d}$.
Thanks to  \cite[Proposition 3.4]{2013Inahama},  it admits that each element $g\in \mathcal{H}^{H,d}$ is
$H'$-H\"older continuous and of finite $(H+1/2)^{-1}<q<2$-variation.

	Then, we 	consider the $\mathbb{R}^{e}$-valued  standard BM $(w_t)_{t\in[0,T]}$,
	$$w_t=(w_t^1,w_t^2,\cdots,w_t^{e}).$$
The reproducing kernel Hilbert space for  $(w_t)_{t\in[0,T]}$, denoted by $\mathcal{H}^{\frac{1}{2},e}$, which is defined as follows, 
	$$\mathcal{H}^{\frac{1}{2},e}:=\big\{ k \in \mathcal{C}_0([0,T],\mathbb{R}^{e})\mid  k _{t}=\int_{0}^{t}  k _{s}^{\prime} d s \text { for  } t\in[0,T] \text { with }\| k \|_{\mathcal{H}^{\frac{1}{2},e}}^{2}:=\int_{0}^{T}| k _{t}^{\prime}|_{\mathbb{R}^{e}}^{2} d t<\infty\big\}.$$

	
In the following, we denote the $\mathbb{R}^{d+e}$-valued	mixed FBM by $(b_t^H, w_t)_{0\le t\le T}$. It is not too difficult to see that $(b^H, w)$  has $H'$-H\"older continuous ($H'\in(0,H)$) and $\lfloor 1 / H \rfloor<{p}<\lfloor 1 / H \rfloor+1$-variation trajectories almost surely. Let $\mathcal{H}:={{\mathcal{H}^{H,d}}\oplus{\mathcal{H}^{\frac{1}{2},e}}}$ be the Cameron-Martin subspace related to $(b_t^H, w_t)_{0\le t\le T}$.  Then, $(\phi, \psi)\in \mathcal{H}$ is of  finite $q$-variation with $(H+1/2)^{-1}<q<2$.

For $N\in\mathbb{N}$, we define 
	\[
	 S_N=\left\{(\phi,\psi) \in \mathcal{H}: 
	\frac{1}{2}\|(\phi,\psi)\|_{\mathcal{H}}^2 := \frac{1}{2} 
	(\|\phi\|_{\mathcal{H}^{H,d}}^{2} +\|\psi \|_{\mathcal{H}^{\frac{1}{2},e}}^{2})
	\leq N\right\}. 
	\]
	The ball $  S_N$ is   a compact Polish space under the weak topology of $\mathcal{H}$. 	
	
	The complete probability space $(\Omega, \mathcal{F}, \mathbb{P})$  supports $b^H$ and $w$ exists independently, where $\Omega=\mathcal{C}_0\left([0,T]: \mathbb{R}^d\right)$,   $\mathbb{P}$ is the unique prabability measure on $\Omega$ and $\mathcal{F}=\mathcal{B}\left(\mathcal{C}_0\left([0,T]: \mathbb{R}^d\right)\right)$ is the $\mathbb{P}$-completion of the Borel $\sigma$-field. Then, we consider the canonical	filtration given by $\left\{\mathcal{F}_t^H: t \in[0,T]\right\}$, where $\mathcal{F}_t^H=\sigma\left\{(b_s^H,w_s): 0 \leq s \leq t\right\} \vee \mathcal{N}$ and $\mathcal{N}$ is the set of the	$\mathbb{P}$-negligible events.

We denote   the set of all 
	$\mathbb{R}^{d}$-valued  processes $(\phi_t,\psi_t)_{t \in [0,T]}$ on 
	$(\Omega, {\mathcal{F}}, {\mathbb{P}})$ by ${\mathcal{A}}_b^N$ for $N\in \mathbb{N}$ and let ${\mathcal{A}}_b =\cup_{N\in \mathbb{N}} {\mathcal{A}}_b^N$.
	Since  each $(\phi,\psi)\in  {\mathcal{A}}_b^N$  is an  random variable taking values in the compact ball ${S}_N$, the family $\{ \mathbb{P}\circ (\phi,\psi)^{-1} :  
	(\phi,\psi) \in {\mathcal{A}}_b^N\}$ of probability measures is  tight automatically.	
	Due to   Girsanov's formula,   for every $(\phi,\psi) \in {\mathcal{A}}_b$,
	the law of  $(b^H +\phi,w+\psi)$ is 
	mutually absolutely continuous to  that of  $(b^H,w)$.
In the following, we recall the  variational representation formula for the mixed FBM, whose precise proof refers to \cite[Proposition 2.3]{2023Inahama}.
\begin{prop}\label{prop2-1}
Let $\alpha\in (0,H)$. For a   bounded Borel measurable function $\Phi:\Omega\to \mathbb{R}$,
\begin{eqnarray}\label{2-1}
-\log \mathbb{E}\left[\exp \left(-\Phi\left(b^H, w\right)\right)\right]=\inf _{(\phi, \psi) \in \mathcal{A}_b} \mathbb{E}\left[\Phi\left(b^H+\phi, w+\psi\right)+\frac{1}{2}\|(\phi, \psi)\|_{\mathcal{H}}^2\right].
\end{eqnarray}
\end{prop} 

		\subsection{Rough Path}\label{sec-2-2}
	In this subsection, we introduce   RP and some explanations which will be utilised in our main proof.  We assume   $\lfloor 1 / H \rfloor<{p}<\lfloor 1 / H \rfloor+1$ and $(H+1/2)^{-1}<q<2$ such that $1/p+1/q>1$. For example, we take $1 / p=H-2\kappa$ and $1 / q=H+1 / 2-\kappa$ with small parameter $0<\kappa<H/2$.

Now, we give the definition of the RP.
	\begin{defn}\cite[Section 2]{2020Friz}\label{RPs}
A continuous map
\begin{eqnarray*}
	\Xi=\big(1, \Xi^{1}, \Xi^{2}\big): \Delta \rightarrow T^{2}(\mathcal{V})=\mathbb{R} \oplus \mathcal{V} \oplus \mathcal{V}^{\otimes 2} ,
\end{eqnarray*}
is said to be a $\mathcal{V}$-valued RP of roughness $2$ if it satisfies the following  conditions,

\textbf{(Condition A)}: For any $s \leq u \leq t$, $\Xi_{s, t}=\Xi_{s, u} \otimes \Xi_{u, t}$ where $\otimes$ stands for the tensor product.

\textbf{(Condition B)}:  $\|\Xi^{1}\|_{\alpha \mathrm{-hld}}<\infty$ and $\|\Xi^{2}\|_{2\alpha \mathrm{-hld}}<\infty$.
	\end{defn}

Obviously, the $0$-th element $1$ is omitted and we denote the RP by $\Xi=\left(\Xi^{1},  \Xi^{2}\right)$. Below, we set $\vertiii{\Xi}_{\alpha-\operatorname{hld}}:=\|\Xi^{1}\|_{\alpha \mathrm{-hld}}+\|\Xi^{2}\|_{2\alpha \mathrm{-hld}}$. 	 The set of all  $\mathcal{V}-$valued RPs with $1/3<\alpha<1/2$ is denoted by $\Omega_{\alpha}(\mathcal{V})$. Equipped with the $\alpha$-H\"older distance, it is a complete space. 	It is easy to verify that  $\Omega_\alpha(\mathcal{V}) \subset \Omega_\beta(\mathcal{V})$ for  $\frac{1}{3}<\beta \leq \alpha \leq \frac{1}{2}$. 
For two different RPs $\Xi=(\Xi^1,\Xi^2)\in \Omega_\alpha(\mathcal{V})$ and $\tilde{\Xi}=(\tilde{\Xi}^1,\tilde{\Xi}^2)\in \Omega_\alpha(\mathcal{V})$, we denote  the  distance between them by $\rho_\alpha(\star, \cdot)$ which is defined as following:
$$\rho_\alpha(\Xi, \tilde{\Xi}):=\|\Xi^1-\tilde{\Xi}^1\|_{\alpha \mathrm{-hld}}+\|\Xi^2-\tilde{\Xi}^2\|_{2\alpha \mathrm{-hld}}.$$


Next, we  give some explanations for RP which will be used in this work.	Firstly, we  show that the  mixed FBM can be lifted to RP, whose precise proof is a minor modification  of  \cite[Proposition 2.2]{2022Yang} by subtracting a term $\frac{1}{2}I_e(t-s)$ where  $I_{e}$ stands  the identity matrix in $\mathbb{R}^{e\times e}$. 
	\begin{rem}\label{rem2}
Let  $(b^H,w)^\mathrm{T}\in \mathbb{R}^{d+e}$ with $H\in(1/3,1/2)$  be the mixed FBM and $\alpha\in (0,H)$.  Then  $(b^H,w)$ can be lifted to RP
$\Lambda=(\Lambda^1,\Lambda^2)\in \Omega_{\alpha}(\mathbb{R}^{d+e})$ with
\begin{eqnarray}\label{2-28}
\Lambda^1_{s t}=\left(b^H_{s t}, w_{s t}\right)^{\mathrm{T}}, \quad \Lambda^2_{s t}=\left(\begin{array}{ll}
B_{s t}^{H,2} & I[b_H, w]_{s t} \\
I[w, b_H]_{s t} & W_{s t}^{2}
\end{array}\right).
\end{eqnarray}
Here,  $(B^{H, 1}, B^{H,2})\in \Omega_{\alpha}(\mathbb{R}^{d})$ is  a canonical geometric RPs associated with FBM  and  $(W^{1}, W^{2})\in \Omega_{\alpha}(\mathbb{R}^{d})$ is a It\^o-type Brownian RP.   Moreover, 
\begin{eqnarray}\label{2-30}
I[b^H, w]_{s t} \triangleq \int_{s}^{t} b^H_{s r} \otimes \mathrm{d}^{\mathrm{I}} w_{r},
\end{eqnarray}
\begin{eqnarray}\label{2-31}
I[w,b^H]_{s t} \triangleq w_{s t} \otimes b^H_{s t}-\int_{s}^{t} \mathrm{~d}^{\mathrm{I}} w_{r} \otimes b^H_{s r},
\end{eqnarray}
where $\int \cdots \mathrm{d}^{\mathrm{I}} w$ stands for the It\^o integral.
	\end{rem}
Moreover, by taking similar way as in \cite[Lemma 4.6]{2022Inahama}, we stress the fact that  $\mathbb{E}\left[\|\Lambda\|_\alpha^q\right]<\infty$ for every $q \in[1, \infty)$.
Then, we turn to the observation that  $u\in \mathcal{H}^{H,d}$   can be lifted to RP.
\begin{rem}\label{rem}
	Let  $H\in (1/3,1/2)$ and $\alpha\in (0,H)$. The 	elements $u\in \mathcal{H}^{H,d}$ can be lifted to RP $U=(U^1,U^2)\in \Omega_{\alpha}(\mathbb{R}^d)$ with
\begin{eqnarray}\label{3-9}
U^1_{s,t}=u_{s,t},\quad U^2_{s,t}=\int_{s}^{t}u_{s,r}du_r
\end{eqnarray}
	where	 $U^2$ is well-defined in the variation setting. Moreover, $U=(U^1,U^2)$ is a locally Lipschitz continuous mapping from $\mathcal{H}^{H,d}$ to $\Omega_\alpha(\mathbb{R}^d)$.
\end{rem}
\para{Proof}. According to the  property that   $u\in \mathcal{H}^{H,d}$  is of finite $(H+1/2)^{-1}<q<2$-variation and Young integral theory, it is not too difficult to derive that  $U^2$ is well-defined in the variation setting. Then,
by applying the fact that $u\in \mathcal{H}^{H,d}$  is	$\alpha$-H\"older continuous, the proof is completed.
\qed

Similarly, we can show that the elements $v\in \mathcal{H}^{\frac{1}{2},e}$ can be lifted to RP $V=(V^1,V^2)\in \Omega_{\alpha}(\mathbb{R}^e)$ with
\begin{eqnarray*}\label{3-10}
	V^1_{s,t}=v_{s,t},\quad V^2_{s,t}=\int_{s}^{t}v_{s,r}dv_r
\end{eqnarray*}
	where	 $U^2$ is well-defined since $v$ is differentiable.
	
Next, we will show that the translation of mixed FBM in the direction $h:=(u,v)\in \mathcal{H}$ can be lifted to RP.
\begin{rem}\label{rem1}
Let   $(b^H+u,w+v)$ be the translation of $(b^H,w)^\mathrm{T}\in \mathbb{R}^{d+e}$ with $H\in(1/3,1/2)$ in the direction $h:=(u,v)\in \mathcal{H}$ and $\alpha\in (0,H)$.	Then, $(b^H+u,w+v)$ can be lifted to RP $T^h(\Lambda)=(T^{h,1}(\Lambda),T^{h,2}(\Lambda))\in \Omega_{\alpha}(\mathbb{R}^{d+e})$, which is defined as following:
		\begin{eqnarray}\label{3-8}
		T_{s,t}^{h,1}(\Lambda)&=&(b^H+u,w+v)_{s,t},	\cr
		T_{s, t}^{h,2}(\Lambda)&=&\left(\begin{array}{ll}
		B^{H,2}+I[b^H,u]+I[u,b^H]+ U^2& I[b_H, w]+I[b_H, v]+I[u, w]+I[u, v] \\
		I[w, b_H]+I[w, u]+	I[v, b_H]+	I[v, u] & W^2+I[w,v]+I[v,w]+V^2
		\end{array}\right)_{s, t}\cr
		&=&\Lambda^2_{s t}+
		\left(\begin{array}{ll}
		I[b^H,u]+I[u,b^H]+ U^2& I[b_H, v]+I[u, w]+I[u, v] \\
	I[w, u]+	I[v, b_H]+	I[v, u] & I[w,v]+I[v,w]+V^2
		\end{array}\right)_{s, t}.
		\end{eqnarray}
		Here, the second term in \eqref{3-8}
is well-defined in the variation setting.
\end{rem}
\para{Proof}. It is obvious that $T^{h,1}(\Lambda)$ is a  translation of mixed FBM in the direction $h:=(u,v)\in \mathcal{H}$ and it is $\alpha$-H\"older continuous. So we  mainly prove the second level path $T^{h,2}(\Lambda)$ is also well-defined. From \remref{rem2} and \remref{rem}, we have shown that $\Lambda^2$, $U^2$ and $V^2$ are well-defined. Hence, we are in the position to show that the remaining terms are also well-defined  in the variation setting. Firstly, we will prove that $I[b^H,u]$ are well-defined in the Young integral. According to \cite[Theorem 2]{2006Friz}, we have that     $b^H$ can be dominated by 
the  function $\omega_1 (s,t):=\|b^H\|^{1/(H-\kappa)}_{(H-\kappa)-hld}{(t-s)}$ for any small $0<\kappa<H$, which is a so-called control function. Similarly, the elements $u$ is dominated by the controll function $\omega_2 (s,t):=\|u\|^{q}_{W^{\delta,p}}{(t-s)}^{\alpha q}$ for $(H+1/2)^{-1}<q<2$ in the sense of \cite[Page 16]{2002Lyons}. The controll function has following super-additivity properties: for $i=1,2$,
\begin{eqnarray}\label{3-11}
\omega_i(s, r)+\omega_i(r, t) \leqslant \omega_i(s, t) \text { with } 0\leq s \leq r \leq t \leq T.
\end{eqnarray}

Let $J_{s,t}=b^H_s(u_t-u_s)$. Then, for $s\leq r\leq t$, we have
\begin{eqnarray*}\label{3-4}
J_{s,r}+J_{r,t}-J_{s,t}&=&b^H_s(u_r-u_s)+b^H_r(u_t-u_r)-b^H_s(u_t-u_s)\cr
&=&(b^H_t-b^H_s)(u_t-u_s).
\end{eqnarray*}
After that, we take a partition $\mathcal{P}=\{s=t_0\leq t_1 \leq ...\leq t_N=t\}$ and denote  
\begin{eqnarray*}\label{3-5}
J_{s,t}(\mathcal{P})=\sum_{i=1}^{N}J_{t_{i-1},t_i},\quad J_{s,t}(\{s,t\})=J_{s,t}.
\end{eqnarray*}
By taking direct computation and using \eqref{3-11}, we obtain
\begin{eqnarray*}\label{3-6}
|J_{s,t}(\mathcal{P})-J_{s,t}(\mathcal{P} \verb|\| \{t_i\})|&\leq&|J_{t_{i-1},t_i}+J_{t_{i},t_{i+1}}-J_{t_{i-1},t_{i+1}}|\cr
&\leq&|(b^H_{t_{i}}-b^H_{t_{i-1}})(u_{t_{i+1}}-u_{t_{i}})|\cr
&\leq& C \{\omega^{^{1/p}}_1 (t_{i-1},t_{i+1})\omega^{^{1/q}}_2 (t_{i-1},t_{i+1})\}\cr
&\leq& {\big(\frac{2}{N}\big)}^{1/p+1/q} \omega^{1/p}_1 (s,t)\omega^{1/q}_2 (s,t).
\end{eqnarray*}
Then, by iterating the above procedure again, we have 
\begin{eqnarray*}\label{3-7}
	|J_{s,t}(\mathcal{P})-J_{s,t}|&\leq& \sum_{k=2}^{N} {\big(\frac{2}{k-1}\big)}^{1/p+1/q} \omega^{1/p}_1 (s,t)\omega^{1/q}_2 (s,t)\cr
	&\leq& 2^{1/p+1/q}\zeta({1/p+1/q})\omega^{1/p}_1 (s,t)\omega^{1/q}_2 (s,t)\cr
	&\leq& 2^{1/p+1/q}\zeta({1/p+1/q})\|b^H\|^{1/p(H-\kappa)}_{(H-\kappa)-hld} \|u\|_{W^{\delta,p}}{(t-s)}^{\alpha+1/p }\cr
		&\leq& C 2^{1/p+1/q}\zeta({1/p+1/q}){(t-s)}^{\alpha+1/p },
\end{eqnarray*}
where $\zeta$ is the Zeta function.
Since $\alpha+1/p >2\alpha$, we verify that second level path  $\int_{s}^{t}b^H_{s,r}du_r$ is of $2\alpha$-H\"older continuous. By the property that the trajectories of $b^H$  is of ${p}$-variation almost surely for $\lfloor 1 / H \rfloor<{p}<\lfloor 1 / H \rfloor+1$ and $u\in \mathcal{H}^{H,d}$ is
 of finite $(H+1/2)^{-1}<q<2$-variation,   $\int_{s}^{t}b^H_{s,r}du_r$ is well-defined in  the Young integral.
Next, by taking similar estimations as above, we can obtain that the other remaining terms are also well-defined in the Young sense.

Moreover, we could verify that  $T^h(\Lambda)=(T^{h,1}(\Lambda),T^{h,2}(\Lambda))$ satisfies  \textbf{(Condition A)} in \defref{RPs} by some direct  computations. Then  we have $T^h(\Lambda)=(T^{h,1}(\Lambda),T^{h,2}(\Lambda))\in \Omega_{\alpha}({\mathbb{R}^{d+e}})$.
The proof is completed. 
\qed

\bigskip
Next, we introduce the  controlled RP. 
Firstly, we recall the definition of controlled RP with respect to the reference RP $\Xi=\left(\Xi^{1},  \Xi^{2}\right)\in \Omega_{\alpha}(\mathcal{V})$. It says that $(Y, Y^{\dagger}, Y^{\sharp})$
is a $\mathcal{W}$-valued controlled RP with respect to $\Xi=\left(\Xi^{1},  \Xi^{2}\right)\in \Omega_{\alpha}(\mathcal{V})$ if it satisfies the following conditions:
$$
Y_t-Y_s=Y_s^{\dagger} \Xi_{s, t}^1+R^{Y}_{s, t}, \quad(s, t) \in \triangle_{[a, b]} 
$$
and 
$$
\left(Y, Y^{\dagger}, R^{Y}\right) \in \mathcal{C}^{\alpha-\operatorname{hld}}([a, b], \mathcal{W}) \times \mathcal{C}^{\alpha-\operatorname{hld}}([a, b], L(\mathcal{V}, \mathcal{W})) \times \mathcal{C}_2^{2\alpha-\operatorname{hld}}([a, b], \mathcal{W}).
$$
Let $\mathcal{Q}_{\Xi}^\alpha([a, b], \mathcal{W})$ stand for the set of all above controlled RPs.  Denote the   semi-norm of controlled RP $(Y, Y^{\dagger}, R^{Y})\in \mathcal{Q}_{\Xi}^\alpha([a, b], \mathcal{W})$  by 
$$\|\left(Y, Y^{\dagger}, R^{Y}\right)\|_{\mathcal{Q}_{\Xi}^\alpha,[a, b]}=\|Y^{\dagger}\|_{{\alpha-\operatorname{hld}},[a, b]}+\|R^{Y}\|_{{2\alpha-\operatorname{hld}},[a, b]}.$$ 
The controlled RP space $\mathcal{Q}_{\Xi}^\alpha([a, b], \mathcal{W})$ is a Banach space equipped with the norm $|Y_a|_{\mathcal{W}}+|Y_a^{\dagger}|_{L(\mathcal{V}, \mathcal{W})}+\|\left(Y, Y^{\dagger}, R^{Y}\right)\|_{\mathcal{Q}_{\Xi}^\alpha,[a, b]}$. In the following, $(Y, Y^{\dagger}, R^{Y})$ is replaced by  $(Y, Y^{\dagger})$ for simplicity.

For two different controlled RPs $(Y, Y^{\dagger})\in \mathcal{Q}_{\Xi}^\alpha([a, b], \mathcal{W})$ and $(\tilde{Y}, \tilde{Y}^{\dagger})\in \mathcal{Q}_{\tilde \Xi}^\alpha([a, b], \mathcal{W})$, we set their distance as follows,
$$d_{\Xi, \tilde{\Xi}, 2 \alpha}\big(Y, Y^{\dagger} ; \tilde{Y}, \tilde{Y}^{\dagger}\big) \stackrel{\text { def }}{=}\big\|Y^{\dagger}-\tilde{Y}^{\dagger}\big\|_{\alpha-\operatorname{hld}}+\big\|R^Y-R^{\tilde{Y}}\big\|_{{2\alpha-\operatorname{hld}}}.$$

In the following, we show that  the integration of controlled RP  against RP is again a controlled RP, whose precise proof refers to \cite[Proposition 3.2]{2022Inahama}.
	\begin{rem}
		Let $1/3<\alpha<1/2$ and $[a,b]\subset[0,T]$. For a RP $\Xi=\left(\Xi^{1},  \Xi^{2}\right)\in \Omega_{\alpha}(\mathcal{V})$ and controlled RP $(Y, Y^{\dagger}) \in \mathcal{Q}_{\Xi}^\alpha([a, b], L(\mathcal{V}, \mathcal{W}))$,  we have 
		$\left(\int_a^{\cdot} Y_u d {\Xi}_u, Y\right) \in \mathcal{Q}_{\Xi}^\alpha([a, b], \mathcal{W})$.
	\end{rem}	
We now turn to the fine property   of the solution to the controlled RDE:
\begin{prop} \label{prop2.5}
Let $\xi\in \mathcal{W}$ and $\Xi=\left(\Xi^{1},  \Xi^{2}\right)\in \Omega_{\alpha}(\mathcal{V})$ with $1/3<\alpha<1/2$. Assume	$(\Psi; \sigma(\Psi))\in \mathcal{Q}_{\Xi}^\beta([0, T],  \mathcal{W})$ with $1/3<\beta<\alpha<1/2$
	be the (unique)  solution to the following RDE
	\begin{eqnarray}\label{2-39}
	d\Psi=f(\Psi_t)dt+\sigma(\Psi_t) d \Xi_t, \quad \Psi_0=\xi \in \mathcal{W}.
	\end{eqnarray}
	Here,  $f$ is globally bounded and Lipschitz continuous function and $\sigma\in C_b^3$. Similarly, let $(\tilde{\Psi}; \sigma(\tilde\Psi)) \in \mathcal{Q}_{\tilde \Xi}^\beta([0, T],  \mathcal{W})$ with initial value $(\tilde \xi, \sigma(\tilde\xi))$. Assume 
	$$\vertiii{\Xi}_{\alpha-\operatorname{hld}},\vertiii{\tilde\Xi}_{\alpha-\operatorname{hld}}\leq M<\infty.$$
	Then, we  have the (local) Lipschitz estimates as following:
	\begin{eqnarray}\label{2-40}
	d_{\Xi, \tilde{\Xi}, 2 \beta}(\Psi, \sigma(\Psi) ; \tilde{\Psi}, \sigma(\tilde{\Psi})) \leq C_M\left(|\xi-\tilde{\xi}|+\rho_\alpha(\Xi, \tilde{\Xi})\right).
	\end{eqnarray}
	and
	\begin{eqnarray}\label{2-41}
	\|\Psi-\tilde{\Psi}\|_{\beta-\operatorname{hld}} \leq C_M\left(|\xi-\tilde{\xi}|+\rho_\alpha(\Xi, \tilde{\Xi})\right).
	\end{eqnarray}
	Here, $C_M=C(M,\alpha,\beta,L_f,\|\sigma\|_{C_b^3})>0$.
\end{prop}
	\para{Proof}. This proposition is a minor modification of \cite[Theorem 8.5]{2020Friz} with the drift term. According to the definition of controlled RP,  we have
	\begin{eqnarray}\label{2-52}
	\|\Psi-\tilde{\Psi}\|_{\beta-\operatorname{hld}} \leq C(d_{\Xi, \tilde{\Xi}, 2 \beta}\big(\Psi, \Psi^{\dagger} ; \tilde{\Psi}, \tilde{\Psi}^{\dagger}\big)+|\xi-\tilde{\xi}|+\rho_\alpha(\Xi, \tilde{\Xi})),
	\end{eqnarray}
	so it only needs to show  \eqref{2-40} and \eqref{2-41} hold.
	
	Let $0<\tau<T$ and we turn to  prove \eqref{2-40} holds in the time interval $[0,\tau]$ firstly.
To this end, we set $\mathcal{M}_{[0, \tau]}^1,\mathcal{M}_{[0, \tau]}^2:\mathcal{Q}_{\Xi}^\beta([0, \tau],  \mathcal{W})\mapsto \mathcal{Q}_{\Xi}^\beta([0, \tau],  \mathcal{W})$ by
	\begin{eqnarray}\label{2-42}
	\mathcal{M}_{[0, \tau]}^1\left(\Psi, \Psi^{\dagger}\right)=\left(\int_0^{\cdot} \sigma\left(\Psi_s\right) d \Xi_s, \sigma(\Psi)\right),\quad\mathcal{M}_{[0, \tau]}^2\left(\Psi, \Psi^{\dagger}\right)=\left(\int_0^{\cdot} f\left(\Psi_s\right) ds, 0\right)
	\end{eqnarray}
	and $(Z, Z^{\dagger}):=\mathcal{M}_{[0, \tau]}^{\xi}:=(\xi, 0)+\mathcal{M}_{[0, \tau]}^1+\mathcal{M}_{[0, \tau]}^2$. Moreover, we stress the fact  that  the  fixed point of $\mathcal{M}_{[0, \tau]}^{\xi}$ is the solution to the \eqref{2-39} on the time interval $[0,\tau]$ for $0<\tau\le T$. Due to the fixed point theorem, we arrive at
	\begin{eqnarray}\label{2-43}
(\Psi, \sigma(\Psi))=\left(\Psi, \Psi^{\dagger}\right)=\left(Z, Z^{\dagger}\right)=(Z, \sigma(\Psi)). 
	\end{eqnarray}
	Abbreviate $\mathcal{I} \Sigma:=Z_{s, t}$ and $\Sigma:=f(\Psi_s)(t-s)+\sigma(\Psi_s)\Xi^1_{s,t}+\sigma^\dagger(\Psi_s)\Xi^2_{s,t}$. Moverover, $\mathcal{I} \tilde\Sigma$ and $\tilde\Sigma$ could be defined  in a similar way with respect to $\tilde{\Psi}$. By some direct computation, we have
	\begin{eqnarray}\label{2-48}
	R_{s, t}^Z&=&Z_{s, t}-Z_s^{\dagger} \Xi_{s, t}\cr
	&=&\int_s^t f (\Psi_r) dr+\int_s^t \sigma (\Psi_r) d\Xi_r-\sigma (\Psi_s) \Xi_{s, t}\cr
	&=&(\mathcal{I} \Sigma)_{s, t}-\Sigma_{s,t}+\sigma^\dagger(\Psi_s)\Xi^2_{s,t}+f(\Psi_s)(t-s).
	\end{eqnarray}
	We set $\mathcal{Q}:=\Sigma-\tilde\Sigma$. After that,  we  obtain that
	\begin{eqnarray}\label{2-49}
	|R_{s, t}^Z-R_{s, t}^{\tilde{Z}}|&=&\big|(\mathcal{I} \mathcal{Q})_{s, t}-\mathcal{Q}_{s, t}\big|+\big|\sigma^\dagger(\Psi_s) \Xi^2_{s, t}-\sigma^\dagger(\tilde\Psi_s) \tilde{\Xi}^2_{s, t}\big|+\big|(f (\Psi_s)-f (\tilde\Psi_s))(t-s)\big|\cr
	&\le& C\|\delta \mathcal{Q}\|_{3 \alpha}|t-s|^{3 \beta}+\big|\sigma^\dagger(\Psi_s) \Xi^2_{s, t}-\sigma^\dagger(\tilde\Psi_s) \tilde{\Xi}^2_{s, t}\big|\cr
	&&+L_f\tau^\beta\|\Psi-\tilde\Psi\|_{\beta-\operatorname{hld}}|t-s|+C|\xi-\tilde{\xi}||t-s|
	\end{eqnarray}
	where $L_f$ is the Lipschitz coefficient of $f$ and $\delta \mathcal{Q}_{s, u, t}=R_{s, u}^{\sigma(\tilde{\Psi})} \tilde{\Xi}^1_{u, t}-R_{s, u}^{\sigma(\Psi)} \Xi^1_{u, t}+\sigma^\dagger(\tilde{\Psi})_{s, u} \tilde{\Xi}^2_{u, t}-\sigma^\dagger(\Psi)_{s, u} \Xi^2_{u, t}$.
	
	Furthermore, a straightforward estimate furnishes that
	\begin{eqnarray}\label{2-53}
|Z_{s,t}^{\dagger}-\tilde Z_{s,t}^{\dagger}|&=&\big|\sigma({Z})_{s, t}-\sigma(\tilde{Z})_{s, t}\big|\cr
&=&\big|\sigma({\Psi})_{s, t}-\sigma(\tilde{\Psi})_{s, t}\big|\cr
&=&\big|(\sigma^\dagger({\Psi})_{0, s}+\sigma^\dagger({\Psi})_{0})\Xi_{s,t}-(\sigma^\dagger(\tilde{\Psi})_{0, s}+\sigma^\dagger(\tilde{\Psi})_{0})\tilde{\Xi}_{s,t}+R_{s, t}^{\sigma({\Psi})}-R_{s, t}^{\sigma(\tilde{\Psi})}\big|\cr
&\le&C|t-s|^\beta\big(|\sigma({\Psi})_{0}-\sigma(\tilde{\Psi})_{0}|+|t-s|^{\alpha-\beta}\|\sigma^\dagger({\Psi})-\sigma^\dagger(\tilde{\Psi})\|_{\beta-\operatorname{hld}}\big.\cr
&&\big.+\rho_\alpha(\Xi, \tilde{\Xi})+\|R^{\sigma({\Psi})}-R^{\sigma(\tilde{\Psi})}\|_{2\beta-\operatorname{hld}}\big).
	\end{eqnarray}
As a consequence of   \cite[Theorem 4.17]{2020Friz} and \eqref{2-49}--\eqref{2-53}, we see that 
	\begin{eqnarray}\label{2-44}
	d_{\Xi, \tilde{\Xi}, 2 \beta}\big(\Psi, \Psi^{\dagger} ; \tilde{\Psi}, \tilde{\Psi}^{\dagger}\big) & =&d_{\Xi, \tilde{\Xi}, 2 \beta}\big(Z, Z^{\dagger} ; \tilde{Z}, \tilde{Z}^{\dagger}\big) \cr
	&=&\|Z^{\dagger}-\tilde Z^{\dagger}\|_{\beta-\operatorname{hld}}+\|R^Z-R^{\tilde{Z}}\|_{2\beta-\operatorname{hld}}\cr
	& \lesssim& \rho_\alpha(\Xi, \tilde{\Xi})+|\xi-\tilde{\xi}|\cr
	&&+\tau^\beta d_{\Xi, \tilde{\Xi}, 2 \beta}\big(\sigma(\Psi), {\sigma^\dagger(\Psi)}
	 ; \sigma(\tilde\Psi), {\sigma^\dagger(\tilde\Psi)}\big)+L_f\tau^{\beta}\|\Psi-\tilde\Psi\|_{\beta-\operatorname{hld}}.
	\end{eqnarray}
	Next, with aid of the \cite[Theorem 7.6]{2020Friz}, we observe that
\begin{eqnarray}\label{2-45}
d_{\Xi, \tilde{\Xi}, 2 \beta}\big(\sigma(\Psi), {\sigma^\dagger(\Psi)}
; \sigma(\tilde\Psi), {\sigma^\dagger(\tilde\Psi)}\big)& \lesssim& \rho_\alpha(\Xi, \tilde{\Xi})+|\xi-\tilde{\xi}|+	d_{\Xi, \tilde{\Xi}, 2 \beta}\big(\Psi, \Psi^{\dagger} ; \tilde{\Psi}, \tilde{\Psi}^{\dagger}\big).
\end{eqnarray}
Therefore, by combining \eqref{2-52} and \eqref{2-44}--\eqref{2-45}, it deduces that  there exists a positive constant $C_M:=C(M,\alpha,\beta,L_f)$ such that 
\begin{eqnarray}\label{2-46}
d_{\Xi, \tilde{\Xi}, 2 \beta}\big(\Psi, \Psi^{\dagger} ; \tilde{\Psi}, \tilde{\Psi}^{\dagger}\big) 
& \le& C_M\big[\rho_\alpha(\Xi, \tilde{\Xi})+|\xi-\tilde{\xi}|+\tau^\beta d_{\Xi, \tilde{\Xi}, 2 \beta}\big(\Psi, \Psi^{\dagger} ; \tilde{\Psi}, \tilde{\Psi}^{\dagger}\big)+\tau^{\beta}\|\Psi-\tilde\Psi\|_{\beta-\operatorname{hld}} \big]\cr
& \le& C_M\big[\rho_\alpha(\Xi, \tilde{\Xi})+|\xi-\tilde{\xi}|+\tau^\beta d_{\Xi, \tilde{\Xi}, 2 \beta}\big(\Psi, \Psi^{\dagger} ; \tilde{\Psi}, \tilde{\Psi}^{\dagger}\big) \big]
\end{eqnarray}
holds.
By taking $\tau>0$ such that $C_M\tau^\beta<1/2$, we find
\begin{eqnarray}\label{2-47}
d_{\Xi, \tilde{\Xi}, 2 \beta}\big(\Psi, \Psi^{\dagger} ; \tilde{\Psi}, \tilde{\Psi}^{\dagger}\big) 
& \le& C_M(\rho_\alpha(\Xi, \tilde{\Xi})+|\xi-\tilde{\xi}|). 
\end{eqnarray}
Then, with   \eqref{2-52},  we arrive at
\begin{eqnarray}\label{2-51}
\|\Psi-\tilde{\Psi}\|_{\beta-\operatorname{hld}} &\leq&C(d_{\Xi, \tilde{\Xi}, 2 \beta}\big(\Psi, \Psi^{\dagger} ; \tilde{\Psi}, \tilde{\Psi}^{\dagger}\big)+|\xi-\tilde{\xi}|+\rho_\alpha(\Xi, \tilde{\Xi}))\cr
&\leq& C_M\big(|\xi-\tilde{\xi}|+\rho_\alpha(\Xi, \tilde{\Xi})\big).
\end{eqnarray}
Taking iteration techniques over $[0,T]$ furnishes that  \eqref{2-40} and \eqref{2-41} hold at the time interval $[0,T]$.
This proof is completed. \qed

		\section{Assumptions and Statement of our Main Result}\label{sec-3}
		In this section, we give necessary assumptions and the statement of our main LDP result. In the all following sections, we set $1/3<\beta<\alpha<H<1/2$.

	We write $Z^{\varepsilon,\delta}=(X^{\varepsilon,\delta},Y^{\varepsilon,\delta})$. Then, the precise definition of slow-fast RDE \eqref{1} can be rewritten as following:
\begin{eqnarray}\label{3-1}
	Z_t^{\varepsilon,\delta}=Z_0+\int_0^t F_{\varepsilon,\delta}\big(Z_s^{\varepsilon,\delta}\big) d s+\int_0^t \Sigma_{\varepsilon,\delta}\big(Z_s^{\varepsilon,\delta}\big) d(\varepsilon\Lambda_s), \quad\big(Z^{\varepsilon,\delta}\big)_t^{\dagger}=\Sigma_{\varepsilon,\delta}\big(Z_t^{\varepsilon,\delta}\big), \quad t \in[0, T].
\end{eqnarray}
with the initial value $Z_0=(X_0,Y_0)$ and
\begin{eqnarray*}\label{3-2}
F_{\varepsilon,\delta}(x, y)=\left(\begin{array}{c}
	f_1(x, y) \\
	\delta^{-1} f_2(x, y)
\end{array}\right), \quad \Sigma_{\varepsilon,\delta}(x, y)=\left(\begin{array}{cc}
	\sigma_1(x) & O \\
	O & (\varepsilon\delta)^{-1 / 2} \sigma_2(x, y)
\end{array}\right) .
\end{eqnarray*}
 Here, $\varepsilon\Lambda=(\sqrt{\varepsilon}\Lambda^1,\varepsilon\Lambda^2)\in \Omega_{\alpha}(\mathbb{R}^{d+e})$ is the dilation of $\Lambda=(\Lambda^1,\Lambda^2)\in \Omega_{\alpha}(\mathbb{R}^{d+e})$, which is defined in \eqref{2-28}.
Then, $(Z^{\varepsilon,\delta},\big(Z^{\varepsilon,\delta}\big)^{\dagger})\in \mathcal{Q}_{\varepsilon\Lambda}^\beta([a, b], \mathbb{R}^{m+n})$ with $1/3<\beta<\alpha<1/2$ is a controlled RP, where the Gubinelli derivative $\big(Z^{\varepsilon,\delta}\big)^{\dagger}$ is defined as following:
\begin{eqnarray*}\label{3-21}
	{(Z^{\varepsilon,\delta})}^\dagger:= \Sigma_{\varepsilon,\delta}(x, y)=\left(\begin{array}{cc}
		\sigma_1(x) & O \\
		O & (\varepsilon\delta)^{-1 / 2} \sigma_2(x, y)
	\end{array}\right) .
\end{eqnarray*}
	To ensure the  existence and uniqueness of solutions to the RDE (\ref{3-1}), we impose the following conditions.
	\begin{itemize}
		\item[\textbf{A1}.] $\sigma_{1}\in \mathcal{C}_b^3$.
		\item[\textbf{A2}.] There exists a constant $L> 0$ such that for  any $ (x_1,y_1) $,  $ (x_2,y_2)\in \mathbb{R}^{m} \times\mathbb{R}^{n}$, 
		\begin{equation*}
		\begin{aligned}
		&\left|f_1\left( x_1, y_1\right)-f_1\left( x_2, y_2\right)\right| +
		\left|f_2\left(x_1, y_1\right)-f_2\left(x_2, y_2\right)\right| \leq L\left(\left|x_1-x_2\right|+\left|y_1-y_2\right|\right), 
		\end{aligned}
		\end{equation*}	
		and
	$$|f_1\left( x_1, y_1\right)|\le L$$	
		hold.
		\item[\textbf{A3}.] There exists a constant $L> 0$ such that for  any $ (x_1,y_1) $,  $ (x_2,y_2)\in \mathbb{R}^{m} \times\mathbb{R}^{n}$, 
		\begin{equation*}
		\begin{aligned}
		&\left|\sigma_2\left(x_1, y_1\right)-\sigma_2\left(x_2, y_2\right)\right|  \leq L\left(\left|x_1-x_2\right|+\left|y_1-y_2\right|\right), 
		\end{aligned}
		\end{equation*}	
		and that, for any $x_1\in \mathbb{R}^m$,
		\begin{equation*}
		\sup_{y_1\in\mathbb{R}^{n}}\left|\sigma_2\left(x_1, y_1\right)\right|  \leq L\left(1+\left|x_1\right|\right)
		\end{equation*}
		hold.
	\end{itemize}		
Under  above   {(\textbf{A1})}--{(\textbf{A3})},  one can deduce from \cite[Remark 3.4]{2022Inahama} that the RDE (\ref{3-1}) has a unique local solution. Thanks to \cite[Proposition 4.8]{2022Inahama}, it deduces  that  the probability that $Z^{\varepsilon, \delta}$ explodes  is zero for all $t\in[0,T]$. Therefore, the RDE (\ref{3-1}) admits a unique  solution $Z^{\varepsilon, \delta}$ globally. 
	
	Furthermore, we have that $\left(X^{\varepsilon,\delta}, \sigma_1\left({X}^{\varepsilon,\delta}\right)\right)\in \mathcal{Q}_{\varepsilon B^H}^\beta\left([0, T], \mathbb{R}^m\right)$  is a unique global solution of the  RDE driven by  $\varepsilon B^H=(\sqrt{\varepsilon}B^{H,1},\varepsilon B^{H,2})$ as following:
	$$
	X_t^{\varepsilon,\delta}=X_0+\int_0^t f_1\big(X_s^{\varepsilon,\delta}, Y_s^{\varepsilon,\delta}\big) d s+\int_0^t \varepsilon\sigma_1\big(X_s^{\varepsilon,\delta}\big) dB^H_s, \quad\big(X^{\varepsilon,\delta}\big)_t^{\dagger}=\sigma_1\big(X_t^{\varepsilon,\delta}\big), \quad t \in[0, T] .
	$$
	Then, by taking similar manner as in \cite[Proposition 4.7]{2022Inahama} and exploiting the conclusion that the probability that $Z^{\varepsilon, \delta}$ explodes  is zero, we  deduce the following conclusion.
For each $0<\delta, \varepsilon\leq 1$, $Y^{\varepsilon,\delta}$ satisfies the It\^o SDE as following:
\begin{eqnarray}
Y_t^{\varepsilon,\delta}=Y_0+\frac{1}{\delta} \int_0^{t} f_2\left(X_s^{\varepsilon,\delta}, Y_s^{\varepsilon,\delta}\right) d s+\frac{1}{\sqrt{\delta}} \int_0^{t} \sigma_2\left(X_s^{\varepsilon,\delta}, Y_s^{\varepsilon,\delta}\right) d^{\mathrm{I}} w_s.
\end{eqnarray}

Then there is a measurable map 
\[
\mathcal{G}^{(\varepsilon, \delta)}: \mathcal{C}_0\left([0,T], \mathbb{R}^d\right) \rightarrow \mathcal{C}^{\beta-\operatorname{hld}}\left([0,T], \mathbb{R}^m\right)
\]
such that
$X^{\varepsilon,\delta}:=\mathcal{G}^{\varepsilon,\delta}(\sqrt \varepsilon b^H, \sqrt \varepsilon w)$. 
Furthermore, to study  an LDP for the slow component in slow-fast RDE \eqref{1}, we assume the following conditions.
\begin{itemize}
	\item[(\textbf{A4}).]  Assume that there exist positive constants $C>0$ and $\beta_i >0 \,(i=1,2)$ such that for any  $ (x,y_1),(x,y_2)\in \mathbb{R}^{m} \times\mathbb{R}^{n}$
	\begin{equation*}
	\begin{aligned}
	2\left\langle y_1-y_2, f_2\left(x, y_1\right)-f_2\left(x, y_2\right)\right\rangle+\left|\sigma_2\left(x, y_1\right)-\sigma_2\left(x, y_2\right)\right|^2 
	&\leq-\beta_1\left|y_1-y_2\right|^2, \\
	2\left\langle y_1, f_2\left(x, y_1\right)\right\rangle+\left|\sigma_2\left(x, y_1\right)\right|^2 & \leq-\beta_2\left|y_1\right|^2+C|x|^2+C
	\end{aligned}
	\end{equation*}
	hold. 
	\end{itemize}	
The Assumptions {(\textbf{A4})}  ensures that the solution to the following It\^o SDE with frozen ${X}$
	\[d\tilde{Y}_t = f_{2}({X}, \tilde{Y}_t) dt +\sigma_2( {X}, \tilde{Y}_t)dw_{t}\]
	possesses a unique invariant probability measure $\mu_{{X}}$, which is  deduced from \cite[Theorem 6.3.2]{DaPratoZabczyk}.

Next, we define  the skeleton equation in the rough sense  as follows
	\begin{eqnarray}\label{3}
d\tilde{X}_t = \bar{f}_1(\tilde{X}_t)dt + \sigma_{1}( \tilde{X}_t)dU_t
\end{eqnarray}	
where $\tilde{X}_0=X_0$, $U=(U^1,U^2)\in \Omega_{\alpha}(\mathbb{R}^d)$  and  $\bar{f}_1(\cdot)=\int_{\mathbb{R}^{n}}f_{1}(\cdot, {\tilde Y})\mu_{\tilde{X}}(d{\tilde Y})$. By taking same  estimates as in \cite[Proposition C.5]{2022Inahama}, we can obtain that $\bar{f}_1$ is   Lipschitz continuous and bounded. Then, it is not too difficult to   see that there exists a unique global solution $(\tilde{X},\tilde{X}^{\dagger})\in \mathcal{Q}_{U}^\beta([0, T], \mathbb{R}^m)$ to the RDE \eqref{3}. Moreover, we have  for $0<\beta<\alpha<H$ that 
$$\|\tilde{X}\|_{\beta-\operatorname{hld}}\le c,$$
with the constant $c>0$  independent of $U$. 
Therefore, we  also define a map
\[
\mathcal{G}^{0}: S_{N} \rightarrow  \mathcal{C}^{\beta-\operatorname{hld}}\left([0,T],\mathbb{R}^m\right)
\]
such that its solution $\tilde{X}=\mathcal{G}^{0}(u, v)$.	

\begin{rem}\label{rem3}
The above RDE \eqref{3} coincides with the Young ODE as following:
\begin{eqnarray}\label{4}
d\tilde{X}_t = \bar{f}_1(\tilde{X}_t)dt + \sigma_{1}( \tilde{X}_t)du_t
\end{eqnarray}	
with $\tilde{X}_t=X_0$ and $\bar{f}_1(\cdot)=\int_{\mathbb{R}^{n}}f_{1}(\cdot, {\tilde Y})\mu_{\tilde{X}}(d{\tilde Y})$. 
For $(H+1/2)^{-1}<q<2$, we have $\|(u,v)\|_{q-\operatorname{var}}<\infty$. According to the Young integral theory, it is easy to verify that  there exists a unique solution $\tilde {X} \in  \mathcal{C}^{p-\operatorname{var}}\left([0,T],\mathbb{R}^d\right)$  to  (\ref{4}) in the Young sense for $(u, v) \in S_{N}$ . Moreover, we have
$$\|\tilde{X}\|_{p-\operatorname{var}}\le c,$$
where the constant $c>0$ is independent of $(u, v)$. 
\end{rem}

Now, we give the statement of our main theorem.	
\begin{thm}\label{thm}
 {Let  $H\in(1/3,1/2)$, fix  $1/3<\beta<H$.	Assume {(\textbf{A1})}--{(\textbf{A4})} and $\delta=o(\varepsilon)$.  Let $\varepsilon \to 0$, the slow component $X^{\varepsilon,\delta}$ of system (\ref{1}) satisfies a LDP on $\mathcal{C}^{\beta-\operatorname{hld}}([0,T],\mathbb{R}^{m})$ with a good rate function $I: \mathcal{C}^{\beta-\operatorname{hld}}([0,T],\mathbb{R}^m)\rightarrow [0, \infty)$}
	\begin{eqnarray*}\label{rate}
		I(\xi) &=& {
			\inf\Big\{\frac{1}{2}\|u\|^2_{\mathcal{H}^{H,d}}~:~{ u\in \mathcal{H}^{H,d} 
				\quad\text{such that} \quad\xi =\mathcal{G}^{0}(u, 0)}\Big\} 
		}
		\cr
		&=& \inf\Big\{\frac{1}{2}\|(u,v)\|^2_{\mathcal{H}} ~:~ {( u,v)\in \mathcal{H}\quad\text{such that} \quad\xi =\mathcal{G}^{0}(u, v)}\Big\},
		\qquad 
		\xi\in \mathcal{C}^{\beta-\operatorname{hld}}\left([0,T], \mathbb{R}^m\right).		
	\end{eqnarray*}
\end{thm}		
	
	\section{Prior estimates}\label{sec-4}
	In this section, we fix $\varepsilon,\delta\in(0,1]$. In the next section, we will let $\varepsilon \to 0$.
	To prove \thmref{thm}, some prior estimates should be given. 
	
Firstly, let $(u^{\varepsilon, \delta}, v^{\varepsilon, \delta})\in \mathcal{A}^{b}$. In order  to apply the variational representation \eqref{2-1}, we give  the following controlled slow-fast RDE  associated to (\ref{1}).
	\begin{eqnarray}\label{2}
	\left
	\{
	\begin{array}{ll}
	d\tilde {X}^{\varepsilon, \delta}_t =& f_{1}(\tilde {X}^{\varepsilon, \delta}_t, \tilde {Y}^{\varepsilon, \delta}_t)dt + \sigma_{1}(\tilde {X}^{\varepsilon, \delta}_t)d[T_{t}^u(\varepsilon B^H) ]\\
	 d\tilde {Y}^{\varepsilon, \delta}_t =&  \frac{1}{\delta}f_{2}( \tilde {X}^{\varepsilon, \delta}_t, \tilde {Y}^{\varepsilon, \delta}_t)dt +\frac{1}{\sqrt{\delta\varepsilon}}\sigma_2(\tilde {X}^{\varepsilon, \delta}_t, \tilde {Y}^{\varepsilon, \delta}_t)dv^{\varepsilon, \delta}_t+ \frac{1}{\sqrt{\delta}}\sigma_2(\tilde {X}^{\varepsilon, \delta}_t, \tilde {Y}^{\varepsilon, \delta}_t)dw_{t}.
	\end{array}
	\right.
	\end{eqnarray}
	Here, $T^u(B^H):=(T^{u,1}(\varepsilon B^H),T^{u,2}(\varepsilon B^H)$ with
			\begin{eqnarray}\label{4-1}
			T_{s,t}^{u,1}(\varepsilon B^H)&=&(\sqrt{\varepsilon}b^H+u^{\varepsilon,\delta})_{s,t}	\cr
			T_{s, t}^{u,2}(\varepsilon B^H)&=&\left(
			\varepsilon B^{H,2}+\sqrt{\varepsilon} I[b^H,u^{\varepsilon,\delta}]+\sqrt{\varepsilon} I[u^{\varepsilon,\delta},b^H]+ U^{\varepsilon,\delta,2}
		\right)_{s, t}.
			\end{eqnarray}
			Here, $(u^{\varepsilon, \delta}, v^{\varepsilon, \delta})\in \mathcal{A}^{b}$ is called a pair of control. 	 It is not too hard to verify that there exists a unique solution $(\tilde {X}^{\varepsilon, \delta},\tilde {Y}^{\varepsilon, \delta})$ to the controlled slow-fast system (\ref{2}). 
			

	For $t\in[0,T]$, we set  $t(\Delta)=\left\lfloor\frac{t}{\Delta}\right\rfloor \Delta$ is the nearest breakpoint preceding $t$.  Then, we construct the auxiliary process as following:
	\begin{eqnarray}\label{4-3}
		 d\hat {Y}^{\varepsilon, \delta}_t &=& \frac{1}{\delta} f_{2}( \tilde {X}^{\varepsilon, \delta}_{t(\Delta)}, \hat {Y}^{\varepsilon, \delta}_t)dt + \frac{1}{{\sqrt \delta }}\sigma_2(\tilde {X}^{\varepsilon, \delta}_{t(\Delta)}, \hat {Y}^{\varepsilon, \delta}_t)dw_{t}.
	\end{eqnarray}
	Now we are in the position to give necessary   prior estimates.
	\begin{lem}\label{lem1}
		Assume (\textbf{A1})--(\textbf{A3}) and let  $\nu\ge 1$ and $N\in\mathbb{N}$. Then, for all  $\varepsilon,\delta\in(0,1]$, we have
		\begin{equation}\label{4-6}
		\mathbb{E}\big[\|\tilde X^{\varepsilon,\delta}\|^\nu_{\beta-\operatorname{hld}}\big] \leq C.
		\end{equation}
		Here, $C$ is a positive constant which depends only on $\nu$ and  $N$.
	\end{lem}
\para{Proof}. The controlled slow component $\tilde X^{\varepsilon,\delta}$
satisfies the following  RDE driven by RP $T^u(\varepsilon B^H)\in \Omega_{\alpha}(\mathbb{R}^{d+e})$ with the initial value $X_0$:
\begin{eqnarray}\label{4-4}
	\tilde {X}^{\varepsilon, \delta}_t =X_0+ \int_{0}^{t}f_{1}(\tilde {X}^{\varepsilon, \delta}_s, \tilde {Y}^{\varepsilon, \delta}_s)ds + \int_{0}^{t}\sigma_{1}(\tilde {X}^{\varepsilon, \delta}_s)d[T_{s}^u(\varepsilon B^H) ],\quad (\tilde {X}^{\varepsilon, \delta}_t)^{\dagger}=\sigma_{1}(\tilde {X}^{\varepsilon, \delta}_t), \quad t\in[0,T].
\end{eqnarray}
For every $(\tilde {X}^{\varepsilon, \delta},(\tilde {X}^{\varepsilon, \delta})^{\dagger})\in  \mathcal{Q}_{T^u(\varepsilon B^H)}^\beta([0, T], \mathbb{R}^{m+n})$, we observe that the right hand side of \eqref{4-4} also belongs to $ \mathcal{Q}_{T^u(\varepsilon B^H)}^\beta([0, T], \mathbb{R}^{m+n})$.
Next, by taking same manner as in \cite[Proposition 3.3]{2022Inahama}, it is clear that
\begin{eqnarray}\label{4-5}
\|\tilde {X}^{\varepsilon, \delta}\|_{\beta-\operatorname{hld}}\le c\{(K+1){(\vertiii{T^u(\varepsilon B^H)}_{\alpha-\operatorname{hld}}+1)}^{\iota}\}
\end{eqnarray}
for constants $c$ and $\iota>0$ which only depends on $\alpha$ and $\beta$. Here, the constant $K:=\|\sigma_1\|_{C_b^3} \vee\|f_1\|_{\infty} \vee L$ where $L$ is defined in (\textbf{A2}). Then, for all $\nu\ge 1$, by taking expectation of $\nu$-moments of \eqref{4-5}, we have
\begin{eqnarray}\label{4-5'}
\mathbb{E}[\|\tilde {X}^{\varepsilon, \delta}\|_{\beta-\operatorname{hld}}^\nu]\le c(K+1)^\nu\mathbb{E}[{(\vertiii{T^u(\varepsilon B^H)}_{\alpha-\operatorname{hld}}+1)}^{\nu\iota}]
\end{eqnarray}
Due to the property that for every $1/3<\alpha<H$ and all $\nu\ge 1$,  $\mathbb{E}[\vertiii{T^u(\varepsilon B^H)}_{\alpha-\operatorname{hld}}^\nu]<\infty$, the estimate \eqref{4-6} is derived.   This proof is completed. 
 \qed
	\begin{lem}\label{lem4}
	Assume (\textbf{A1})--(\textbf{A4}) and let $N\in\mathbb{N}$. Then, for every $(u^{\varepsilon},v^{\varepsilon})\in \mathcal{A}_b^N$, we have
	\begin{eqnarray}\label{4.1}
	\int_{0}^{T} \mathbb{E}\big[|\tilde  Y_t^{\varepsilon,\delta}|^2\big] dt\leq C.
	\end{eqnarray}
	Here, $C$ is a positive constant which depends only  on $N$.
\end{lem}
\para{Proof}.
Due to that $Y^{\varepsilon,\delta}$ satisfies the It\^o SDE and by using   It\^o's formula, we have
\begin{eqnarray}\label{4.1'}
{|\tilde  Y_t^{\varepsilon,\delta} |^2} &=& {| Y_0 |^2} + \frac{2}{\delta }\int_0^t {\langle \tilde Y_s^{\varepsilon,\delta},f_2( \tilde {X}_{s}^{\varepsilon,\delta},\tilde Y_s^{\varepsilon,\delta})\rangle ds} 
+ \frac{2}{\sqrt{\delta} }\int_0^t {\langle \tilde {Y}_{s}^{\varepsilon,\delta},\sigma_{2}( {\tilde {X}_{s}^{\varepsilon,\delta},\tilde {Y}_{s}^{\varepsilon,\delta}} ) dw_s\rangle} \cr
&&+ \frac{2}{{\sqrt { \varepsilon\delta } }}\int_0^t {\langle \tilde {Y}_{s}^{\varepsilon,\delta},\sigma_{2} ( {\tilde {X}_{s}^{\varepsilon,\delta},\tilde {Y}_{s}^{\varepsilon,\delta}} ){\frac{dv_s^{\varepsilon,\delta}}{ds}  }\rangle ds}   + \frac{1}{\delta }\int_0^t |{\sigma_{2} }( {\tilde {X}_{s}^{\varepsilon,\delta},\tilde {Y}_{s}^{\varepsilon,\delta}} )|^2ds.
\end{eqnarray}
By similar estimates in  \cite[Lemma 4.3]{2023Inahama}, we derive  that 
for every $(u^{\varepsilon,\delta},v^{\varepsilon,\delta})\in \mathcal{A}_b^N$, 
$\sup_{0\le s\le t}|\tilde  Y_s^{\varepsilon,\delta}|$ has moments of all orders. 
Then,  with the aid of  \lemref{lem1}, we can prove that  the third term in right hand side of (\ref{4.1'}) is a true martingale and  $\mathbb{E}[\int_0^t {\langle \tilde {Y}_{s}^{\varepsilon,\delta},\sigma_{2}( {\tilde {X}_{s}^{\varepsilon,\delta},\tilde {Y}_{s}^{\varepsilon,\delta}} ) dW_s\rangle}]=0$.
Taking expectation for  \eqref{4.1'}, we have
\begin{eqnarray}\label{3.21}
\begin{aligned}
\frac{d\mathbb{E}[{| \tilde {Y}_{t}^{\varepsilon} |^2}]}{dt} = \frac{2}{\delta }\mathbb{E} \big[{\langle \tilde {Y}_{t}^{\varepsilon,\delta},f_2( \tilde {X}_{t}^{\varepsilon,\delta},\tilde {Y}_{t}^{\varepsilon,\delta} )\rangle } \big]
+ \frac{2}{{\sqrt { \varepsilon \delta} }}\mathbb{E} \big[{\langle \tilde {Y}_{t}^{\varepsilon,\delta},\sigma_{2} ( {\tilde {Y}_{t}^{\varepsilon,\delta},\tilde {Y}_{t}^{\varepsilon,\delta}} ){\frac{dv_t^{\varepsilon,\delta}}{dt}  }\rangle }\big]   
+ \frac{1}{\delta }\mathbb{E} \big[|{\sigma_{2} }( {\tilde {X}_{t}^{\varepsilon,\delta},\tilde {Y}_{t}^{\varepsilon,\delta}} )|^2\big].
\end{aligned}
\end{eqnarray} 
By  (\textbf{A4}), we arrive at
\begin{eqnarray}\label{lemma3.12}
\frac{2}{\delta }{\langle \tilde {Y}_{t}^{\varepsilon,\delta},f_2( \tilde {X}_{t}^{\varepsilon,\delta},\tilde {Y}_{t}^{\varepsilon,\delta} )\rangle } +\frac{1}{\delta } |{\sigma_{2} }( {\tilde {X}_{t}^{\varepsilon,\delta},\tilde {Y}_{t}^{\varepsilon,\delta}} )|^2\le- \frac{{ \beta_2 }}{\delta }{{| \tilde {Y}_{t}^{\varepsilon,\delta} |^2}}   + \frac{{C}}{\delta }|\tilde {X}_{t}^{\varepsilon,\delta}|^2+ \frac{{C}}{\delta }.	
\end{eqnarray}
With aid of    (\textbf{A4}) and \lemref{lem1}, we obtain
\begin{eqnarray}\label{3.22}
\frac{2}{{\sqrt { \varepsilon\delta } }}{\langle \tilde {Y}_{t}^{\varepsilon,\delta},\sigma_{2} ( {\tilde {X}_{t}^{\varepsilon,\delta},\tilde {Y}_{t}^{\varepsilon,\delta}} ){\frac{dv_t^{\varepsilon,\delta}}{dt}  }\rangle }&\le&   \frac{L}{{\sqrt { \varepsilon\delta } }}\big( 1 + | \tilde {X}_{t}^{\varepsilon,\delta} |^2 \big)| {\frac{dv_t^{\varepsilon,\delta}}{dt}  } |^2+   \frac{1}{{\sqrt { \varepsilon\delta } }}| \tilde {Y}_{t}^{\varepsilon,\delta} |^2 \cr
&\le&\frac{L}{{\sqrt { \varepsilon\delta } }}\big( 1 + T^2\| \tilde {X}^{\varepsilon,\delta} \|_{\beta-\operatorname{hld}}^2 \big)| {\frac{dv_t^{\varepsilon,\delta}}{dt}  } |^2+   \frac{1}{{\sqrt { \varepsilon\delta } }}| \tilde {Y}_{t}^{\varepsilon,\delta} |^2
\end{eqnarray}
Thus, combine (\ref{3.21})--(\ref{3.22}), it deduces that
\begin{eqnarray*}\label{3.23}
	\begin{aligned}
		\frac{d\mathbb{E}[{| \tilde {Y}_{t}^{\varepsilon,\delta} |^2}]}{dt} 
		&\le  {\frac{{ - \beta_2 }}{2\delta } } \mathbb{E}[{| \tilde {Y}_{t}^{\varepsilon,\delta} |^2}]  +   \frac{{{L T^2}}}{{\sqrt { \varepsilon\delta } }}  \mathbb{E}[{{\| \tilde {X}^{\varepsilon,\delta} \|_{\beta-\operatorname{hld}}^2{| {\frac{dv_t^{\varepsilon,\delta}}{dt}  } |^2}}}]  +   \frac{{{L }}}{\sqrt { \varepsilon\delta } } \mathbb{E}[{| {\frac{dv_t^{\varepsilon,\delta}}{dt}  } |^2}] +\frac{{C}}{\delta }\mathbb{E}[|\tilde {X}_{t}^{\varepsilon,\delta}|^2]+ \frac{{ C}}{\delta }.
	\end{aligned}
\end{eqnarray*}
Furthermore,  by applying the comparison theorem for all $t$, we get  
\begin{eqnarray}\label{3.231}
\begin{aligned}
\mathbb{E}[{| \tilde {Y}_{t}^{\varepsilon,\delta} |^2}] &\le|Y_0|^2 e^{-\frac{\beta_2}{2\delta} t} +   \frac{{{L T^2}}}{{\sqrt { \varepsilon\delta } }}\int_{0}^{t}  e^{-\frac{\beta_2}{2\delta} (t-s)}{\mathbb{E}[{\| \tilde {X}^{\varepsilon,\delta} \|_{\beta-\operatorname{hld}}^2{| {\frac{dv_s^{\varepsilon,\delta}}{ds}  } |^2}}]} ds  \\
&\qquad+   \frac{{{L }}}{\sqrt { \varepsilon\delta } }\int_{0}^{t}  e^{-\frac{\beta_2}{2\delta} (t-s)} { \mathbb{E}[|{\frac{dv_s^{\varepsilon,\delta}}{ds}  } |^2]}ds +\frac{{ C}}{\delta }{\mathbb{E}[{\| \tilde {X}^{\varepsilon,\delta} \|_{\beta-\operatorname{hld}}^2}]}\int_{0}^{t}  e^{-\frac{\beta_2}{2\delta} (t-s)}ds\\
&\qquad+ \frac{{ C}}{\delta }\int_{0}^{t}  e^{-\frac{\beta_2}{2\delta} (t-s)}ds.
\end{aligned}
\end{eqnarray}
Next, by integrating  of (\ref{3.231}) and using the Fubini theorem and \lemref{lem1}, we  can prove that 
\begin{eqnarray*}\label{3.232}
	\begin{aligned}
		\int_{0}^{T}\mathbb{E}[{| \tilde {Y}_{t}^{\varepsilon} |^2}]dt 
		&\le |Y_0|^2 \int_{0}^{T} e^{-\frac{\beta_2}{2\delta} t} dt+   \frac{{{LT^2}}}{{\sqrt { \varepsilon \delta} }}\int_{0}^{T}\int_{0}^{t}  e^{-\frac{\beta_2}{2\varepsilon} (t-s)}{\mathbb{E}[{\| \tilde {X}^{\varepsilon,\delta} \|_{\beta-\operatorname{hld}}^2{| {\frac{dv_s^{\varepsilon,\delta}}{ds}  } |^2}}]} dsdt  \\
		&\quad+   \frac{{{L }}}{\sqrt {\delta \varepsilon } }\int_{0}^{T}\int_{0}^{t}  e^{-\frac{\beta_2}{2\varepsilon} (t-s)} { \mathbb{E}[|{\frac{dv_s^{\varepsilon,\delta}}{ds}  } |^2]}ds + \frac{{ C}}{\delta }\int_{0}^{t}  e^{-\frac{\beta_2}{2\varepsilon} (t-s)}dsdt\\
		&\le |Y_0|^2 e^{-\frac{\beta_2}{2\delta} T} +   \frac{{{LT^2  }}}{{\sqrt { \varepsilon\delta } }}\mathbb{E}\big[\| \tilde {X}^{\varepsilon,\delta} \|_{\beta-\operatorname{hld}}^2\times| \int_{0}^{T}\int_{s}^{T}  e^{-\frac{\beta_2}{2\delta} (t-s)} dt| {\frac{dv_s^{\varepsilon,\delta}}{ds}  } |^2ds\big]  \\
		&\quad+   \frac{{{L }}}{\sqrt { \varepsilon\delta } }\int_{0}^{T}\int_{s}^{T}  e^{-\frac{\beta_2}{2\delta} (t-s)}dt { \mathbb{E}[|{\frac{dv_s^{\varepsilon\delta}}{ds}  } |^2]}ds + \frac{{ C}}{\delta }\int_{0}^{t}  e^{-\frac{\beta_2}{2\delta} (t-s)}ds\\
		&\le |Y_0|^2 e^{-\frac{\beta_2}{2\delta} T} +   \frac{2LT^2 \sqrt{\delta}}{\beta_2\sqrt { \varepsilon } }\mathbb{E}\big[\| \tilde {X}^{\varepsilon,\delta} \|_{\beta-\operatorname{hld}}^2\times| \int_{0}^{T}  e^{-\frac{\beta_2}{2\delta} (T-s)} | {\frac{dv_s^{\varepsilon,\delta}}{ds}  } |^2ds\big]  \\
		&\quad+   \frac{2 L \sqrt{\delta}}{\beta_2\sqrt { \varepsilon }}\int_{0}^{T} e^{-\frac{\beta_2}{2\delta} (T-s)} { \mathbb{E}[|{\frac{dv_s^{\varepsilon\delta}}{ds}  } |^2]}ds +{C}\mathbb{E}[\| \tilde {X}^{\varepsilon,\delta} \|_{\beta-\operatorname{hld}}^2]\int_{0}^{T} e^{-\frac{\beta_2}{2\delta} (T-s)}ds+ C.
	\end{aligned}
\end{eqnarray*}
By using the condition  that   $0<\delta<\varepsilon
\le 1$ and $(u^{\varepsilon}, v^{\varepsilon})\in {\mathcal{A}}_{b}^N$, we derive
\begin{eqnarray*}\label{3.24}
	\int_{0}^{T}\mathbb{E}[{| \tilde {Y}_{t}^{\varepsilon,\delta} |^2}]dt &\le& {C}\mathbb{E}[\| \tilde {X}^{\varepsilon,\delta} \|_{\beta-\operatorname{hld}}^2]+C.
\end{eqnarray*}
Thus, by exploiting the  \lemref{lem1}, the estimate (\ref{4.1}) follows at once. The proof is completed.
\qed 
\begin{lem}\label{lem3}
Assume (\textbf{A1})--(\textbf{A4}) and let $N\in\mathbb{N}$, we have 
	$$\mathbb{E}\big[\big|\tilde Y_{t}^{\varepsilon,\delta}-\hat Y_{t}^{\varepsilon,\delta}\big|^2\big] \le C(\frac{\sqrt{\delta}}{{\sqrt{\varepsilon}}}+\Delta^{2\beta}).$$
	Here, $C>0$ is a  constant which depends only  on $N,\alpha,\beta$.
\end{lem}
\para{Proof}. By It\^o's formula, we have
\begin{eqnarray}\label{4-7}
\mathbb{E}[{| \tilde Y_{t}^{\varepsilon,\delta}-\hat Y_{t}^{\varepsilon,\delta}|^2}] &=&  \frac{2}{\delta }\mathbb{E}\bigg[\int_0^t {\langle\tilde Y_{s}^{\varepsilon,\delta}-\hat Y_{s}^{\varepsilon,\delta},f_2( \tilde {X}_{s}^{\varepsilon,\delta},\tilde Y_s^{\varepsilon,\delta})-f_2( \tilde {X}_{s(\Delta)}^{\varepsilon,\delta},\hat Y_s^{\varepsilon,\delta})\rangle ds}\bigg]\cr
&&+ \frac{1}{\delta }\mathbb{E}\bigg[\int_0^t |{\sigma_{2} }( {\tilde {X}_{s}^{\varepsilon,\delta},\tilde {Y}_{s}^{\varepsilon,\delta}} )-{\sigma_{2} }( {\tilde {X}_{s(\Delta)}^{\varepsilon,\delta},\hat {Y}_{s}^{\varepsilon,\delta}} )|^2ds\bigg] \cr
&&+ \frac{2}{{\sqrt { \varepsilon\delta } }}\mathbb{E}\bigg[\int_0^t {\langle \tilde Y_{s}^{\varepsilon,\delta}-\hat Y_{s}^{\varepsilon,\delta},\sigma_{2} ( {\tilde {X}_{s}^{\varepsilon,\delta},\tilde {Y}_{s}^{\varepsilon,\delta}} ){\frac{dv_s^{\varepsilon,\delta}}{ds}  }\rangle ds}   \bigg].
\end{eqnarray}
By differentiating with respect to $t$ for  \eqref{4-7}, we find that
\begin{eqnarray}\label{4-8}
\frac{d}{dt}\mathbb{E}[{| \tilde Y_{t}^{\varepsilon,\delta}-\hat Y_{t}^{\varepsilon,\delta}|^2}] &=&  \frac{2}{\delta }\mathbb{E}\big[ {\langle\tilde Y_{t}^{\varepsilon,\delta}-\hat Y_{t}^{\varepsilon,\delta},f_2( \tilde {X}_{t}^{\varepsilon,\delta},\tilde Y_t^{\varepsilon,\delta})-f_2( \tilde {X}_{t(\Delta)}^{\varepsilon,\delta},\hat Y_t^{\varepsilon,\delta})\rangle }\big]\cr
&&+ \frac{1}{\delta }\mathbb{E}\big[ |{\sigma_{2} }( {\tilde {X}_{t}^{\varepsilon,\delta},\tilde {Y}_{t}^{\varepsilon,\delta}} )-{\sigma_{2} }( {\tilde {X}_{t(\Delta)}^{\varepsilon,\delta},\hat {Y}_{t}^{\varepsilon,\delta}} )|^2\big] \cr
&&+ \frac{2}{{\sqrt { \varepsilon\delta } }}\mathbb{E}\big[ {\langle \tilde Y_{t}^{\varepsilon,\delta}-\hat Y_{t}^{\varepsilon,\delta},\sigma_{2} ( {\tilde {X}_{t(\Delta)}^{\varepsilon,\delta},\tilde {Y}_{t}^{\varepsilon,\delta}} ){\frac{dv_s^{\varepsilon,\delta}}{dt}  }\rangle }   \big]\cr
&=&  \frac{1}{\delta }\mathbb{E}\big[ {2\langle\tilde Y_{t}^{\varepsilon,\delta}-\hat Y_{t}^{\varepsilon,\delta},f_2( \tilde {X}_{t}^{\varepsilon,\delta},\tilde Y_t^{\varepsilon,\delta})-f_2( \tilde {X}_{t}^{\varepsilon,\delta},\hat Y_t^{\varepsilon,\delta})\rangle }+|{\sigma_{2} }( {\tilde {X}_{t}^{\varepsilon,\delta},\tilde {Y}_{t}^{\varepsilon,\delta}} )-{\sigma_{2} }( {\tilde {X}_{t}^{\varepsilon,\delta},\hat {Y}_{t}^{\varepsilon,\delta}} )|^2\big]\cr
&& + \frac{2}{\delta }\mathbb{E}\big[ {\langle\tilde Y_{t}^{\varepsilon,\delta}-\hat Y_{t}^{\varepsilon,\delta},f_2( \tilde {X}_{t}^{\varepsilon,\delta},\hat Y_t^{\varepsilon,\delta})-f_2( \hat {X}_{t(\Delta)}^{\varepsilon,\delta},\hat Y_t^{\varepsilon,\delta})\rangle }\big]\cr
&&+ \frac{2}{\delta }\mathbb{E}\big[ \langle{\sigma_{2} }( {\tilde {X}_{t}^{\varepsilon,\delta},\tilde {Y}_{t}^{\varepsilon,\delta}} )-{\sigma_{2} }( {\tilde {X}_{t}^{\varepsilon,\delta},\hat {Y}_{t}^{\varepsilon,\delta}} ),{\sigma_{2} }( {\tilde {X}_{t}^{\varepsilon,\delta},\hat {Y}_{t}^{\varepsilon,\delta}} )-{\sigma_{2} }( {\tilde {X}_{t(\Delta)}^{\varepsilon,\delta},\hat {Y}_{t}^{\varepsilon,\delta}} )\rangle\big] \cr
&&+ \frac{1}{\delta }\mathbb{E}\big[ |{\sigma_{2} }( {\tilde {X}_{t}^{\varepsilon,\delta},\hat {Y}_{t}^{\varepsilon,\delta}} )-{\sigma_{2} }( {\tilde {X}_{t(\Delta)}^{\varepsilon,\delta},\hat {Y}_{t}^{\varepsilon,\delta}} )|^2\big] \cr
&&+ \frac{2}{{\sqrt { \varepsilon\delta } }}\mathbb{E}\big[ {\langle \tilde Y_{t}^{\varepsilon,\delta}-\hat Y_{t}^{\varepsilon,\delta},\sigma_{2} ( {\tilde {X}_{t}^{\varepsilon,\delta},\tilde {Y}_{t}^{\varepsilon,\delta}} ){\frac{dv_t^{\varepsilon,\delta}}{dt}  }\rangle }   \big]\cr
&=:&I_1+I_2+I_3+I_4+I_5.
\end{eqnarray}
For the first term $I_1$, by using (\textbf{A4}), we obtian that
\begin{eqnarray}\label{4-9}
I_1 \le -\frac{\beta_1}{\delta}\mathbb{E}\big[\big|\tilde Y_{t}^{\varepsilon,\delta}-\hat Y_{t}^{\varepsilon,\delta}\big|^2\big].
\end{eqnarray}
Then, we compute the second term $I_2$ by using (\textbf{A2}) and \lemref{lem1} as follows,
\begin{eqnarray}\label{4-10}
I_2 &\le& \frac{C_1}{\delta}\mathbb{E}\big[\big|\tilde Y_{t}^{\varepsilon,\delta}-\hat Y_{t}^{\varepsilon,\delta}\big|\cdot\big|\tilde X_{t}^{\varepsilon,\delta}-\tilde X_{t(\Delta)}^{\varepsilon,\delta}\big|\big]\cr
&\le&\frac{\beta_1}{4\delta}\mathbb{E}\big[\big|\tilde Y_{t}^{\varepsilon,\delta}-\hat Y_{t}^{\varepsilon,\delta}\big|^2\big]+\frac{C_2}{\delta}\mathbb{E}\big[\big|\tilde X_{t}^{\varepsilon,\delta}-\tilde X_{t(\Delta)}^{\varepsilon,\delta}\big|^2\big]\cr
&\le& \frac{\beta_1}{4\delta}\mathbb{E}\big[\big|\tilde Y_{t}^{\varepsilon,\delta}-\hat Y_{t}^{\varepsilon,\delta}\big|^2\big]+\frac{C_2}{\delta}\Delta^{2\beta}\mathbb{E}[\| \tilde {X}^{\varepsilon,\delta} \|_{\beta-\operatorname{hld}}^2].
\end{eqnarray}
where $C_1,C_2>0$ is independent of $\varepsilon,\delta$.

For the third term $I_3$ and forth term $I_4$, we estimate them as following:
\begin{eqnarray}\label{4-11}
I_3+I_4 &\le& \frac{C}{\delta}\mathbb{E}\big[\big|\tilde Y_{t}^{\varepsilon,\delta}-\hat Y_{t}^{\varepsilon,\delta}\big|\cdot\big|\tilde X_{t}^{\varepsilon,\delta}-\tilde X_{t(\Delta)}^{\varepsilon,\delta}\big|+\big|\tilde X_{t}^{\varepsilon,\delta}-\tilde X_{t(\Delta)}^{\varepsilon,\delta}\big|^2\big]\cr
&\le& \frac{\beta_1}{4\delta}\mathbb{E}\big[\big|\tilde Y_{t}^{\varepsilon,\delta}-\hat Y_{t}^{\varepsilon,\delta}\big|^2\big]+\frac{C_3}{\delta}\mathbb{E}\big[\big|\tilde X_{t}^{\varepsilon,\delta}-\tilde X_{t(\Delta)}^{\varepsilon,\delta}\big|^2\big]\cr
&\le& \frac{\beta_1}{4\delta}\mathbb{E}\big[\big|\tilde Y_{t}^{\varepsilon,\delta}-\hat Y_{t}^{\varepsilon,\delta}\big|^2\big]+\frac{C_3}{\delta}\Delta^{2\beta}\mathbb{E}[\| \tilde {X}^{\varepsilon,\delta} \|_{\beta-\operatorname{hld}}^2],
\end{eqnarray}
where $C_3>0$ is independent of $\varepsilon,\delta$. Here, for the first inequality, we used  (\textbf{A3}). For the final inequality, we applied \lemref{lem1} and the definition of H\"older norm. 

For the fifth term $I_5$, by applying (\textbf{A3}), we  derive
\begin{eqnarray}\label{4-12}
I_5 &\le& \frac{C}{\sqrt{\varepsilon\delta}}\mathbb{E}\big[\big|\tilde Y_{t}^{\varepsilon,\delta}-\hat Y_{t}^{\varepsilon,\delta}\big|\times \big|1+\tilde X_{t}^{\varepsilon,\delta}\big|\big|{\frac{dv_t^{\varepsilon,\delta}}{dt}  }\big|\big]\cr
&\le& \frac{\beta_1}{4\sqrt{\varepsilon\delta}}\mathbb{E}\big[\big|\tilde Y_{t}^{\varepsilon,\delta}-\hat Y_{t}^{\varepsilon,\delta}\big|^2\big]+\frac{C_4}{\sqrt{\varepsilon\delta}}\mathbb{E}\big[\big|1+\tilde X_{t}^{\varepsilon,\delta}\big|^2\big|{\frac{dv_t^{\varepsilon,\delta}}{dt}  }\big|^2\big],
\end{eqnarray}
where $C_4>0$ is independent of $\varepsilon,\delta$.
Then, by combining \eqref{4-8}--\eqref{4-12}, we have
\begin{eqnarray}\label{4-13}
\frac{d}{dt}\mathbb{E}[{| \tilde Y_{t}^{\varepsilon,\delta}-\hat Y_{t}^{\varepsilon,\delta}|^2}]&\le&-\frac{\beta_1}{4\delta}\mathbb{E}\big[\big|\tilde Y_{t}^{\varepsilon,\delta}-\hat Y_{t}^{\varepsilon,\delta}\big|^2\big]+\frac{C_4}{\sqrt{\varepsilon\delta}}\mathbb{E}\big[\big|1+\tilde X_{t}^{\varepsilon,\delta}\big|^2\big|{\frac{dv_t^{\varepsilon,\delta}}{dt}  }\big|^2\big]+\frac{C_2+C_3}{\delta}\Delta^{2\beta}.
\end{eqnarray}
Thanks to  the Gronwall inequality \cite[Lemma A.1 (2)]{2022Inahama} and \lemref{lem1}, we can observe that 
\begin{eqnarray}\label{4-14}
\mathbb{E}[{| \tilde Y_{t}^{\varepsilon,\delta}-\hat Y_{t}^{\varepsilon,\delta}|^2}]&\le&\frac{C_4\sqrt{\delta}}{\sqrt{\varepsilon}}\int_{0}^{t}\mathbb{E}\big[\big|1+\tilde X_{t}^{\varepsilon,\delta}\big|^2\big|{\frac{dv_t^{\varepsilon,\delta}}{dt}  }\big|^2dt\big]+(C_2+C_3)\Delta^{2\beta}T\cr
&\le&\frac{C_5\sqrt{\delta}}{\sqrt{\varepsilon}}\mathbb{E}\big(1+\| \tilde {X}^{\varepsilon,\delta} \|_{\beta-\operatorname{hld}}^{2\beta} T^{2\beta}\big)+(C_2+C_3)\Delta^{2\beta}T\cr
&\le&C(\frac{\sqrt{\delta}}{{\sqrt{\varepsilon}}}+\Delta^{2\beta}).
\end{eqnarray}
The proof is completed.
\qed

\section{Proof of \thmref{thm}}
In this section, we are ultimately going to  prove  our main result \thmref{thm}. We divide this proof into three steps.

	\textbf{Step 1}.
	The proof is deterministic in this step .
	Let $(u^{(j)}, v^{(j)}),(u, v)\in S_N$  such that $(u^{(j)}, v^{(j)})\rightarrow(u, v)$ as $j\rightarrow\infty$ with the weak topology in $\mathcal{H}$. 
	In this step,  we will  prove that  
	\begin{eqnarray} \label{3.26}
	\mathcal{G}^{0}( u^{(j)} ,v^{(j)} )\rightarrow\mathcal{G}^{0}(u,v) \label{step2}
	\end{eqnarray}
	in $\mathcal{C}^{\beta-\operatorname{hld}}([0,T],\mathbb{R}^{m})$ as $j \to \infty$.
	
The skeleton equation satisfies the RDE as follows
	\begin{eqnarray}\label{5-1}
	d\tilde{X}^{(j)}_t = \bar{f}_1(\tilde{X}^{(j)}_t)dt + \sigma_{1}( \tilde{X}^{(j)}_t)dU^{(j)}_t
	\end{eqnarray}	
	where $\tilde{X}^{(j)}_t=X_0$, $U^{(j)}=({(U^{(j)})^1},{(U^{(j)})^2})\in \Omega_{\alpha}(\mathbb{R}^d)$  and  $\bar{f}_1(\cdot)=\int_{\mathbb{R}^{n}}f_{1}(\cdot, {\tilde Y})\mu_{\cdot}(d{\tilde Y})$. By the conclusion that $\bar{f}_1$ is   Lipschitz continuous and bounded and using \cite[Proposition 3.3]{2022Inahama}, we   obtain that there exists a unique global solution $(\tilde{X}^{(j)},(\tilde{X}^{(j)})^{\dagger})\in \mathcal{Q}_{U}^\beta([0, T], \mathbb{R}^m)$ to the \eqref{5-1}. Moreover, we have
	$$\|\tilde{X}^{(j)}\|_{\beta-\operatorname{hld}}\le c$$
	holds for $0<\beta<\alpha<H$. Here, the constant $c>0$ which is  independent of $U$. 
	
	Due to  a compact embedding  $\mathcal{C}^{\beta-\operatorname{hld}}([0,T],\mathbb{R}^{m}) \subset  \mathcal{C}^{(\beta-\theta)-\operatorname{hld}}([0,T],\mathbb{R}^{m})$  for any small parameter $0<\theta<\beta$, we have that the family $\{\tilde{X}^{(j)}\}_{ j\ge 1}$ is pre-compact in $ \mathcal{C}^{(\beta-\theta)-\operatorname{hld}}([0,T],\mathbb{R}^{m})$.  Let $\tilde X$ be any limit point. Then, there exists  a subsequence of $\{\tilde{X}^{(j)}\}_{ j\ge 1}$ (denoted by the same symbol) weakly converges to  $\tilde X$ in $ \mathcal{C}^{(\beta-\theta)-\operatorname{hld}}([0,T],\mathbb{R}^{m})$. In the following, 
		we will prove   that the limit point $\tilde X$  satisfies the  RDE as follows,
		\begin{eqnarray}\label{5-2}
	d\tilde{X}_t = \bar{f}_1(\tilde{X}_t)dt + \sigma_{1}( \tilde{X}_t)dU_t.
	\end{eqnarray}
	
According to  \remref{rem3}, we could emphasize that $\{\tilde{X}^{(j)}\}_{ j\ge 1}$ coincides with the   following   ODE:
	\begin{eqnarray}\label{5-3}
	d\tilde{X}^{(j)}_t = \bar{f}_1(\tilde{X}^{(j)}_t)dt + \sigma_{1}( \tilde{X}^{(j)}_t)du^{(j)}_t.
	\end{eqnarray}	
where $\|(u^{(j)},v^{(j)})\|_{q-\operatorname{var}}<\infty$ with $(H+1/2)^{-1}<q<2$ for all $j\ge 1$. 
Due to the Young integral theory, it is not too difficult to verify that  for all  $(u, v) \in S_{N}$, there exists a unique solution $\{\tilde{X}^{(j)}\}_{ j\ge 1} \in  \mathcal{C}^{p-\operatorname{var}}\left([0,T],\mathbb{R}^m\right)$  to  \eqref{5-3} in the Young sense. In fact, $\{\tilde{X}^{(j)}\}_{ j\ge 1}$ is independent of $\{v^{(j)}\}_{ j\ge 1}$. Moreover, we have
$$\|\tilde{X}^{(j)}\|_{p-\operatorname{var}}\le c,$$
where the constant $c>0$ is independent of $(u^{(j)}, v^{(j)})$. 
Note that the Young integral  $u^{(j)}\mapsto\int_{0}^{\cdot}\sigma_{1}(\tilde{X}_s^{(j)})du^{(j)}_s$ is a linear continuous map from $\mathcal{H}^d$ to $\mathcal{C}^{p-\operatorname{var}}\left([0,T],\mathbb{R}^m\right)$.	
	
Let us show that the limit point $\tilde{X}$ satisfies the skeleton equation \eqref{4}. By the direct computation, we derive
\begin{eqnarray}\label{5-4}
\big|\tilde{X}^{(j)}_t-\tilde{X}_t\big|&\le& \bigg|\int_{0}^{t} \big[\bar{f}_1(\tilde{X}^{(j)}_s)-\bar{f}_1(\tilde{X}_s)\big]ds\bigg|+\bigg|\int_{0}^{t} \big[\sigma_{1}( \tilde{X}^{(j)}_s)-\sigma_{1}( \tilde{X}_s)\big]du^{(j)}_s\bigg|\cr
&&+\bigg|\int_{0}^{t} \sigma_{1}( \tilde{X}_s)\big[du^{(j)}_s-du_s\big]\bigg|\cr
&=:& J_1+J_2+J_3.
\end{eqnarray}
For the first term $J_1$, by using the result that $\bar{f}_1$ is   Lipschitz continuous and bounded, we have
\begin{eqnarray}\label{5-5}
J_1\le L\int_{0}^{t} |\tilde{X}^{(j)}_s-\tilde{X}_s|ds\le C\sup_{0\le s\le t}|\tilde{X}^{(j)}_s-\tilde{X}_s|.
\end{eqnarray}
After that, by applying (\textbf{A1}), we estimate $J_2$ as following:
\begin{eqnarray}\label{5-6}
J_2\le C\bigg|\int_{0}^{t} |\tilde{X}^{(j)}_s-\tilde{X}_s|du^{(j)}_s\bigg|\le CT\|u^{(j)}\|_{q-\operatorname{var}}\sup_{0 \leq t \leq T}|\tilde{X}^{(j)}_t-\tilde{X}_t|\le C_1\sup_{0 \leq t \leq T}|\tilde{X}^{(j)}_t-\tilde{X}_t|
\end{eqnarray}
where $C_1>0$ only depends on $N$ and $q$.
Since  $\{\tilde X^{(j)}\}_{ j\ge 1}$ converges to $\tilde X$ in the uniform norm, it is an immediate consequence that
$J_1+ J_2\to 0$ as $j \to \infty$.

Next, it proceeds to estimates $J_3$. To do this, we set $B(u^{(j)},\tilde{X}):=\int_{0}^{t} \sigma_{1}( \tilde{X}_s)du^{(j)}_s$, which is a bilinear continuous map from $\mathcal{H}^{H,d}\times \mathcal{C}^{p-\operatorname{var}}\left([0,T],\mathbb{R}^m\right)$ to $\mathbb{R}$.   According to the  Riesz representation theorem,  there exists a unique element in $  \mathcal{H}^{H,d}$ (denoted by $B(\cdot,\tilde{X})$) such that  $B(u^{(j)},\tilde{X})=\langle B(\cdot,\tilde{X}),u^{(j)}\rangle_{\mathcal{H}^{H,d}}$ for all $u^{(j)}\in \mathcal{H}^{H,d}$. Note that $B(\cdot,\tilde{X})\in (\mathcal{H}^{H,d})^{*}\cong \mathcal{H}^{H,d}$. Then, we have 
\begin{eqnarray}\label{5-7}
J_3&=& |B(u^{(j)},\tilde{X})-B(u,\tilde{X})|\cr
&=& |\langle B(\cdot,\tilde{X}),u^{(j)}\rangle_{\mathcal{H}^{H,d}}-\langle B(\cdot,\tilde{X}),u\rangle_{\mathcal{H}^{H,d}}|.
\end{eqnarray}
Since $(u^{(j)}, v^{(j)})\rightarrow(u, v)$ as $j\rightarrow\infty$ with the weak topology in $\mathcal{H}$, we prove that  $J_3$ converges to 0 as $j\rightarrow\infty$.

By combining \eqref{5-4}--\eqref{5-7} and \remref{rem3}, it is clear  that the limit point $\tilde X$ satisfies the ODE (\ref{4}). Consequently, we obtain that  $\{\tilde X^{(j)}\}_{ j\ge 1}$ weakly converges to $\tilde X$ in $\mathcal{C}^{(\beta-\theta)-\operatorname{hld}}([0,T],\mathbb{R}^{m})$  for any small $0<\theta<\beta$.

\textbf{Step 2}. 
We carry out probabilistic arguments in this step. 
Let $0<N<\infty$ and 
assume $0<\delta=o(\varepsilon)\le 1$ and we will take $\varepsilon \to 0$.

Assume  $(u^{\varepsilon, \delta}, v^{\varepsilon, \delta})\in \mathcal{A}^{N}_b$ such that $(u^{\varepsilon, \delta}, v^{\varepsilon, \delta})$ weakly converges  to $(u, v)$ as $\varepsilon \to 0$. 
In this step, we will prove that   $\tilde {X}^{\varepsilon, \delta}$ weakly converges to $\tilde {X}$ in $\mathcal{C}^{\beta-\operatorname{hld}}([0,T],\mathbb{R}^{m})$ as $\varepsilon\rightarrow 0$, that is,
\begin{eqnarray} \label{step3}
\mathcal{G}^{(\varepsilon,\delta)}(\sqrt \varepsilon b^H+u^{\varepsilon, \delta}, \sqrt \varepsilon w+v^{\varepsilon, \delta})\xrightarrow{\textrm{weakly}}\mathcal{G}^{0}(u , v) \quad \text{as }\varepsilon\rightarrow 0.
\end{eqnarray}	
We rewrite  the controlled slow component of RDE (\ref{2}) as following,
\[
\tilde {X}^{\varepsilon,\delta}:=\mathcal{G}^{(\varepsilon,\delta)}(\sqrt \varepsilon b^H+u^{\varepsilon,\delta}, \sqrt \varepsilon w+v^{\varepsilon,\delta}).
\]
Before showing \eqref{step3} hold, we define an auxiliary process $\hat  X^{\varepsilon,\delta}$  satisfies the following RDE:
\begin{eqnarray}\label{5-8}
d\hat {X}^{\varepsilon,\delta}_t = \bar f_{1}(\tilde {X}^{\varepsilon,\delta}_t)dt +\sigma_{1}(\tilde {X}^{\varepsilon,\delta}_t)d[T_{t}^u(\varepsilon B^H) ]
\end{eqnarray}
with initial value $\hat {X}^{\varepsilon,\delta}_0=X_0$.
By taking similar manner as in \lemref{lem1}, we can have
\begin{eqnarray}\label{5-35}
\mathbb{E}[\|\hat {X}^{\varepsilon,\delta}\|^2_{\beta-\operatorname{hld}}]\le   C
\end{eqnarray}
where $C>0$ only depends on $\alpha,\beta$ and $N$.

Now, we are in the position to  give some prior estimates which will be used in proving \eqref{step3}. Firstly, by some direct computation, we can get that
\begin{eqnarray}\label{5-9}
\tilde {X}^{\varepsilon,\delta}_t-\hat {X}^{\varepsilon,\delta}_t&=&\int_{0}^{t}[f_{1}(\tilde {X}^{\varepsilon, \delta}_s, \tilde {Y}^{\varepsilon, \delta}_s)-f_{1}(\tilde {X}^{\varepsilon, \delta}_{s(\Delta)}, \tilde {Y}^{\varepsilon, \delta}_s)]dt +\int_{0}^{t}[f_{1}(\tilde {X}^{\varepsilon, \delta}_{s(\Delta)}, \tilde {Y}^{\varepsilon, \delta}_s)-f_{1}(\tilde {X}^{\varepsilon, \delta}_{s(\Delta)}, \hat {Y}^{\varepsilon, \delta}_s)]ds \cr
&&+\int_{0}^{t}[f_{1}(\tilde {X}^{\varepsilon, \delta}_{s(\Delta)}, \hat {Y}^{\varepsilon, \delta}_s)-\bar f_{1}(\tilde {X}^{\varepsilon, \delta}_{s(\Delta)})]ds +\int_{0}^{t}[\bar f_{1}(\tilde {X}^{\varepsilon, \delta}_{s(\Delta)})-\bar f_{1}(\tilde {X}^{\varepsilon, \delta}_s)]ds \cr
&& +\int_{0}^{t}[\bar f_{1}(\tilde {X}^{\varepsilon,\delta}_{s})-\bar f_{1}(\hat {X}^{\varepsilon,\delta}_{s})]ds+\int_{0}^{t}[\sigma_{1}(\tilde {X}^{\varepsilon, \delta}_s)-\sigma_{1}(\hat {X}^{\varepsilon,\delta}_s)]d[T_{s}^u(\varepsilon B^H) ]\cr
&:=&K_1+K_2+K_3+K_4+K_5+K_6.
\end{eqnarray}
Firstly, we estimate $K_1$ with H\"older inequality, (\textbf{A2}) and \lemref{lem1}, 
\begin{eqnarray}\label{5-11}
\mathbb{E}[\sup_{0 \leq t \leq T}|K_1|^2] &=&\mathbb{E}\bigg[\sup_{0 \leq t \leq T}\big|\int_{0}^{t}[f_{1}(\tilde {X}^{\varepsilon, \delta}_s, \tilde {Y}^{\varepsilon, \delta}_s)-f_{1}(\tilde {X}^{\varepsilon, \delta}_{s(\Delta)}, \tilde {Y}^{\varepsilon, \delta}_s)]ds\big|^2\bigg]\cr
&\le&LT\int_{0}^{T}\mathbb{E}[|\tilde {X}^{\varepsilon, \delta}_s-\tilde {X}^{\varepsilon, \delta}_{s(\Delta)}|^2]ds\cr
&\le&LT^2 \mathbb{E}[\|\tilde {X}^{\varepsilon, \delta}\|^2_{\beta-\operatorname{hld}}]\Delta^{2\beta}.
\end{eqnarray}
For the second term $K_2$, with aid of the H\"older inequality and  \lemref{lem3}, we get 
\begin{eqnarray}\label{5-12}
\mathbb{E}[\sup_{0 \leq t \leq T}|K_2|^2] &=&\mathbb{E}\bigg[\sup_{0 \leq t \leq T}\big|\int_{0}^{t}[f_{1}(\tilde {X}^{\varepsilon, \delta}_{s(\Delta)}, \tilde {Y}^{\varepsilon, \delta}_s)-f_{1}(\tilde {X}^{\varepsilon, \delta}_{s(\Delta)}, \hat {Y}^{\varepsilon, \delta}_s)]ds\big|^2\bigg] \cr
&\le&TL	\int_{0}^{T}\mathbb{E}\big[\big|\tilde Y_{s}^{\varepsilon,\delta}-\hat Y_{s}^{\varepsilon,\delta}\big|^2\big]ds \le C(\frac{\sqrt{\delta}}{{\sqrt{\varepsilon}}}+\Delta^{2\beta}).
\end{eqnarray}
In the following part, we will estimate $K_3$. To this end, we set $M_{s,t}=\int_{s}^{t}[f_{1}(\tilde {X}^{\varepsilon, \delta}_{r(\Delta)}, \hat {Y}^{\varepsilon, \delta}_r)-\bar f_{1}(\tilde {X}^{\varepsilon, \delta}_{r(\Delta)})]dr$. Then, we give some prior estimates. Set $1/2<\eta<1$.
When $0<t-s<2\Delta$, it is immediate to see that
\begin{eqnarray}\label{5-18}
|M_{s,t}|&\le&L {(2\Delta)^{1-\eta}}{(t-s)^\eta}
\end{eqnarray}
When $t-s>2\Delta$, by using the Schwarz inequality, we obtain
\begin{eqnarray}\label{5-19}
	\frac{\left|M_{s, t}\right|^2}{(t-s)^{2\eta}} & \le&\frac{\big|M_{s,(\lfloor s / \Delta\rfloor+1) \Delta}+\sum_{k=\lfloor s / \Delta\rfloor+1}^{\lfloor t / \Delta\rfloor-1} M_{k \Delta,(k+1) \Delta}+M_{\lfloor t / \Delta\rfloor \Delta, t}\big|^2}{(t-s)^{2\eta}}  \cr
	& \le&  C \Delta^{2-2\eta}+\frac{2C(t-s)^{1-2\eta}}{\Delta} \sum_{k=0}^{\lfloor T / \Delta\rfloor-1}\left|M_{k \Delta,(k+1) \Delta}\right|^2.
\end{eqnarray}
Then, by \eqref{5-18} and \eqref{5-19}, it deduces that
\begin{eqnarray}\label{5-20}
\mathbb{E}[\|K_3\|_{\beta-\operatorname{hld}}^2]&=&\mathbb{E}\bigg[\bigg\|\int_{0}^{\cdot}[f_{1}(\tilde {X}^{\varepsilon, \delta}_{s(\Delta)}, \hat {Y}^{\varepsilon, \delta}_s)-\bar f_{1}(\tilde {X}^{\varepsilon, \delta}_{s(\Delta)})]ds\bigg\|_{\beta-\operatorname{hld}}^2\bigg]\cr
&\le &C\Delta^{2(1-\eta)}+\frac{CT}{\Delta^{(1+2\eta)}} \max _{0 \leq k \leq\left\lfloor\frac{T}{\Delta}\right\rfloor-1} \mathbb{E}\big[\big|\int_{k \Delta}^{(k+1) \Delta}\big(f_1( \tilde X_{k \Delta}^{\varepsilon, \delta},\hat {Y}^{\varepsilon, \delta}_s)-\bar{f}_1( \tilde X_{k \Delta}^{\varepsilon, \delta})\big) d s\big|^2\big]
\end{eqnarray}
According to some direct but cumbersome computation, we arrive at
\begin{eqnarray}\label{5-13}
&&\max _{0 \leq k \leq\left\lfloor\frac{T}{\Delta}\right\rfloor-1} \mathbb{E}\big[\big|\int_{k \Delta}^{(k+1) \Delta}\big(f_1( \tilde X_{k \Delta}^{\varepsilon, \delta},\hat {Y}^{\varepsilon, \delta}_s)-\bar{f}_1( \tilde X_{k \Delta}^{\varepsilon, \delta})\big) d s\big|^2\big]\cr
&\leq & C \delta^2 \max _{0 \leq k \leq\left\lfloor\frac{T}{\Delta}\right\rfloor-1} \int_0^{\frac{\Delta}{\delta}} \int_r^{\frac{\Delta}{\delta}} \mathbb{E}\big[\big\langle f_1( \tilde X_{k \Delta}^{\varepsilon,\delta}, \hat{Y}_{s \varepsilon+k \Delta}^{(\varepsilon,\delta)})-\bar{f}_1( \tilde X_{k \Delta}^{\varepsilon,\delta}),f_1( \tilde X_{k \Delta}^{\varepsilon,\delta}, \hat{Y}_{r \varepsilon+k \Delta}^{\varepsilon,\delta})-\bar{f}_1( \tilde X_{k \Delta}^{\varepsilon,\delta})\big\rangle\big] d s d r\cr
& \leq& C \delta^2 \max _{0 \leq k \leq\left\lfloor\frac{T}{\Delta}\right\rfloor-1} \int_0^{\frac{\Delta}{\varepsilon}} \int_r^{\frac{\Delta}{\varepsilon}} e^{-\frac{\beta_1}{2}(s-r)} d s d r \cr
& \leq &C \delta^2\left(\frac{2}{\beta_1} \frac{\Delta}{\delta}-\frac{4}{\beta_1^2}+e^{\frac{-\beta_1}{2} \frac{\Delta}{\delta}}\right)\cr
& \leq &C \delta\Delta.
\end{eqnarray}
Here, we exploit the exponential ergodicity of $\hat Y^{\varepsilon,\delta}$, that is \begin{eqnarray}\label{3.56}
\mathbb{E}\big[\big\langle f_1( \tilde x_{k \Delta}^{\varepsilon}, \hat{y}_{s \varepsilon+k \Delta}^{\varepsilon})-\bar{f}_1( \tilde x_{k \Delta}^{\varepsilon}),f_1( \tilde x_{k \Delta}^{\varepsilon}, \hat{y}_{r \varepsilon+k \Delta}^{\varepsilon})-\bar{f}_1( \tilde x_{k \Delta}^{\varepsilon})\big\rangle\big]\leq C e^{-\frac{\beta_1}{2}(s-\zeta)},
\end{eqnarray}
where $\beta_1$ is in (\textbf{A4}), whose precise proof refers to  \cite[Appendix B]{2023Pei} for instance.
So we have
\begin{eqnarray}\label{5-21}
\mathbb{E}[\|K_3\|_{\beta-\operatorname{hld}}^2]\le C\Delta^{2(1-\eta)}+\frac{CT\delta}{\Delta^{2\eta}}.
\end{eqnarray}
Next, for the forth term $K_4$, by applying  that $\bar{f}_1$ is   Lipschitz continuous and bounded, we obtain
\begin{eqnarray}\label{5-14}
\mathbb{E}[\sup_{0 \leq t \leq T}|K_4|^2] &=&\mathbb{E}\bigg[\sup_{0 \leq t \leq T}\bigg|\int_{0}^{t}[\bar f_{1}(\tilde {X}^{\varepsilon, \delta}_{s(\Delta)})-\bar f_{1}(\tilde {X}^{\varepsilon, \delta}_s)]ds\bigg|^2\bigg]\cr
&\le& LT\int_{0}^{T}\mathbb{E}[|\tilde {X}^{\varepsilon, \delta}_{s(\Delta)}-\tilde {X}^{\varepsilon, \delta}_{s}|^2]ds\cr
&\le&LT^2 \mathbb{E}[\|\tilde {X}^{\varepsilon, \delta}\|^2_{\beta-\operatorname{hld}}]\Delta^{2\beta}.
\end{eqnarray}
Next, we set
\begin{eqnarray}\label{5-30}
Q_t&:=&(\tilde {X}^{\varepsilon,\delta}_t-\hat {X}^{\varepsilon,\delta}_t)-\bigg\{\int_{0}^{t}[\bar f_{1}(\tilde {X}^{\varepsilon, \delta}_s)-\bar f_{1}(\hat {X}^{\varepsilon,\delta}_s)]ds\bigg\} \cr
&&-\bigg\{\int_{0}^{t}[\sigma_{1}(\tilde {X}^{\varepsilon, \delta}_s)-\sigma_{1}(\hat {X}^{\varepsilon,\delta}_s)]d[T_{s}^u(\varepsilon B^H)]\bigg\}
\end{eqnarray}
The estimates \eqref{5-11}--\eqref{5-30} furnish the following observation that $Q\in \mathcal{C}^{1-\operatorname{hld}}([0,T],\mathbb{R}^m)$ and
\begin{eqnarray}\label{5-32}
\mathbb{E}\left[\|Q\|_{2 \beta}^2\right] \leq C\left(\Delta^{2 \beta}+\Delta^{2(1-2 \beta)}+\Delta^{-4 \beta} \delta+\frac{\sqrt{\delta}}{\sqrt{\varepsilon}}\right) .
\end{eqnarray}
Due to  \cite[Proposition 3.5]{2022Inahama}, it deduces that there exist positive constants $c$ and $\nu$ such that
\begin{eqnarray}\label{5-31}
\|\tilde {X}^{\varepsilon, \delta}-\hat {X}^{\varepsilon,\delta}\|_{\beta-\operatorname{hld}} \leq c \exp \left[c\left(K^{\prime}+1\right)^\nu\big(\vertiii{T^u(\varepsilon B^H)}_{\alpha-\operatorname{hld}}+1\big)^\nu\right]\|Q\|_{2\beta-\operatorname{hld}}.
\end{eqnarray}
Here, $K^{\prime}=\max \{\|\sigma_1\|_{C_b^3},\|f_1\|_{\infty}, L\}$.
Then, we choose some suitable $\Delta>0$ such that $\mathbb{E}[\|Q\|_{2\beta-\operatorname{hld}}^2] \to 0$ as $\varepsilon$. For instance, we could choose $\Delta:=\delta^{1 /(4 \beta)} \log \delta^{-1}$.
Therefore, we have that $\|\tilde {X}^{\varepsilon, \delta}-\hat {X}^{\varepsilon,\delta}\|^2_{\beta-\operatorname{hld}}$ converges to $0$ in probability as $\varepsilon \to 0$. 

On the other hand, with \lemref{lem1} and \eqref{5-35}, it is clear to find that
\begin{eqnarray}\label{5-33}
\mathbb{E}[\|\tilde {X}^{\varepsilon, \delta}-\hat {X}^{\varepsilon,\delta}\|^2_{\beta-\operatorname{hld}}] \leq c \mathbb{E}[\|\tilde {X}^{\varepsilon, \delta}\|^2_{\beta-\operatorname{hld}}]+\mathbb{E}[\|\hat {X}^{\varepsilon,\delta}\|^2_{\beta-\operatorname{hld}}]\leq C.
\end{eqnarray}
So it shows that  $\|\tilde {X}^{\varepsilon, \delta}-\hat {X}^{\varepsilon,\delta}\|^2_{\beta-\operatorname{hld}}$ is uniformly integrable.
Then, as a consequence of  the bounded convergence theorem, we prove that  $\mathbb{E}[\|\tilde {X}^{\varepsilon, \delta}-\hat {X}^{\varepsilon,\delta}\|^2_{\beta-\operatorname{hld}}]$ converges to $0$  as $\varepsilon \to 0$. 

Then, we define 
\begin{eqnarray}\label{5-54}
d\tilde {X}^{\varepsilon}_t = \bar f_{1}(\tilde {X}^{\varepsilon}_t)dt +\sigma_{1}(\tilde {X}^{\varepsilon}_t)dU^{\varepsilon,\delta}_t 
\end{eqnarray}
with initial value $\tilde {X}^{\varepsilon}_0=X_0$.
By taking similar manner as in \lemref{lem1}, we observe
\begin{eqnarray}\label{5-53}
\mathbb{E}[\|\tilde {X}^{\varepsilon}\|^2_{\beta-\operatorname{hld}}]\le   C
\end{eqnarray}
where $C>0$ only depends on $\alpha,\beta$ and $N$.

By using \propref{prop2.5}, we have that 
\begin{eqnarray}\label{5-40}
\|\hat {X}^{\varepsilon, \delta}-\tilde {X}^{\varepsilon}\|_{\beta-\operatorname{hld}} &\le& C_{N,B^H}\rho_\alpha(T^u(\varepsilon B^H), U^{\varepsilon,\delta})\cr
&\le&C_{N,B^H}(\|\sqrt{\varepsilon}b^H\|_{\alpha-\operatorname{hld}}+\|\varepsilon I[b^H,u^{\varepsilon,\delta}] \|_{2\alpha-\operatorname{hld}})\cr
&&+C_{N,B^H}(\|\varepsilon I[u^{\varepsilon,\delta},b^H] \|_{2\alpha-\operatorname{hld}}+\|\varepsilon B^{H,2}\|_{2\alpha-\operatorname{hld}})\cr
&\le&C_{N,B^H}\sqrt{\varepsilon}
\end{eqnarray}
where $C_{N,B^H}:=C_{N,\vertiii{B^{H}}_{\alpha-\operatorname{hld}}}>0$ is independent of $\varepsilon$ and $\delta$. On the other hand,  it is not too intractable to verify that
\begin{eqnarray}\label{5-41}
\mathbb{E}[\|\hat {X}^{\varepsilon, \delta}-\tilde {X}^{\varepsilon}\|^2_{\beta-\operatorname{hld}}] \leq 2 \mathbb{E}[\|\hat {X}^{\varepsilon, \delta}\|^2_{\beta-\operatorname{hld}}]+2\mathbb{E}[\|\tilde {X}^{\varepsilon}\|^2_{\beta-\operatorname{hld}}]\leq C.
\end{eqnarray}
So it implies that  $\|\hat {X}^{\varepsilon, \delta}-\tilde {X}^{\varepsilon}\|^2_{\beta-\operatorname{hld}}$ is uniformly integrable.
Then, by applying the bounded convergence theorem, it is immediate to see that  $\mathbb{E}[\|\hat {X}^{\varepsilon, \delta}-\tilde {X}^{\varepsilon}\|^2_{\beta-\operatorname{hld}}]$ converges to $0$  as $\varepsilon \to 0$. 

In the following, we will show  that $\tilde X^\varepsilon$ converges in distribution to $\tilde X$ as $\varepsilon \to 0$. By  \remref{rem} and condition that $(u^{\varepsilon, \delta}, v^{\varepsilon, \delta})\in  \mathcal{A}^{N}_b$, we have that $ U^{\varepsilon,\delta}: \mathcal{H}^{H,d} \mapsto \Omega_\alpha(\mathbb{R}^d)$ is a Lipschitz continuous mapping. Next, by \propref{prop2.5}, we obtain that  $\tilde {X}^{\varepsilon}$ is a continuous solution map with respect to RP $ U^{\varepsilon,\delta}$. With aid of the condition that  $(u^{\varepsilon, \delta}, v^{\varepsilon, \delta})$ weakly converges  to $(u, v)$ as $\varepsilon \to 0$ and continuous mapping theorem, it deduces that $\tilde X^\varepsilon$ converges in distribution to $\tilde X$ as $\varepsilon \to 0$.

By employing the   Portemanteau theorem \cite[Theorem 13.16]{2020Klenke}, we have for any bounded Lipschitz functions $F:\mathcal{C}^{\beta-\operatorname{hld}}\left([0,T], \mathbb{R}^m\right) \to\mathbb{R}$,  that  
	\begin{eqnarray*}\label{5-34}
		|\mathbb{E}[F(\tilde X^{\varepsilon,\delta})]-\mathbb{E}[F(\tilde X)]|&\le& |\mathbb{E}[F(\tilde X^{\varepsilon,\delta})]-\mathbb{E}[F(\tilde X^{\varepsilon})]|+|\mathbb{E}[F(\tilde X^{\varepsilon})]-\mathbb{E}[F(\tilde X)]|\cr
		&\le& \|F\|_\textrm{Lip}\mathbb{E}[\|\tilde X^{\varepsilon,\delta}-\tilde X^{\varepsilon}\|_{\beta\textrm{-hld}}^2]^{\frac{1}{2}}+|\mathbb{E}[F(\tilde X^{\varepsilon})]-\mathbb{E}[F(\tilde X)]|\to 0
	\end{eqnarray*} 
	as $\varepsilon \to 0$. Here, $\|F\|_\textrm{Lip}$ is the Lipschitz constant of $F$. So we have  proved (\ref{step3}).

\textbf{Step 3}. 
By \textbf{Step 1} and \textbf{Step 2},   for   every bounded and continuous     function $\Phi: \mathcal{C}^{\beta-\operatorname{hld}}([0,T],\mathbb{R}^m)\to \mathbb{R}$, we have that the Laplace   lower bound 
	\begin{eqnarray}\label{5-36}
	\liminf _{\varepsilon \rightarrow 0} -\varepsilon \log \mathbb{E}\big[e^{-\frac{\Phi(X^{\varepsilon,\delta})}{\varepsilon} }\big] \geq\inf _{\psi:=\mathcal{G}^0(u,v) \in \mathcal{C}^{{\beta-\operatorname{hld}}}([0,T],\mathbb{R}^m)}\left[\Phi(\psi)+I(\psi)\right]
	\end{eqnarray}
	and the Laplace upper bound
	\begin{eqnarray}\label{5-37}
	\limsup _{\varepsilon \rightarrow 0} -\varepsilon\log \mathbb{E}\big[e^{-\frac{\Phi(X^{\varepsilon,\delta})}{\varepsilon} }\big] \leq\inf _{\psi:=\mathcal{G}^0(u,v)  \in \mathcal{C}^{{\beta-\operatorname{hld}}}([0,T],\mathbb{R}^m)}\left[\Phi(\psi)+I(\psi)\right]
	\end{eqnarray}
	hold and the goodness of rate function $I$. The precise proof for \eqref{5-36}--\eqref{5-37} refers to  \cite[Theorem 3.1]{2023Inahama} as an example.

Hence, our LDP result is concluded by the equivalence between the LDP and Laplace priciple at once. 
This proof is completed.\qed

	\section*{Acknowledgments}
     This work was partly supported by the Key International (Regional) Cooperative Research Projects of the NSF of China (Grant 12120101002), the NSF of China (Grant 12072264), the Fundamental Research Funds for the Central Universities, the Research Funds for Interdisciplinary Subject of Northwestern Polytechnical University,  the Shaanxi Provincial Key R\&D Program (Grants 2020KW-013, 2019TD-010).

\end{document}